\documentclass[letterpaper,oneside,10pt]{article}

\usepackage[USenglish]{babel} 
\usepackage[T1]{fontenc}
\usepackage[ansinew]{inputenc}

\usepackage{lmodern} 

\usepackage{graphicx} 

\usepackage{float}
\usepackage{caption}
\usepackage{amsmath}
\usepackage{amsthm}
\usepackage{amsfonts}
\usepackage{amssymb}
\usepackage{mathtools}
\usepackage[hyphens]{url}

\let\Right\right 
\let\Left\left 
\makeatletter 
\def\right#1{\Right#1\@ifnextchar){\!\right}{}} 
\def\left#1{\Left#1\@ifnextchar({\!\left}{}} 
\makeatother

\begin{document}

\pagestyle{empty} 


\title{On Lerch's formula for the Fermat quotient}
\author{John Blythe Dobson (j.dobson@uwinnipeg.ca)}

\maketitle

\tableofcontents 
\cleardoublepage 

\pagestyle{plain} 


\begin{abstract}
\noindent
This paper explores some consequences of Lerch's 1905 formula for the Fermat quotient, with special attention to the sums $s(k, N) = \sum_{j = \lfloor kp/N\rfloor + 1}^{\lfloor(k + 1)p/N\rfloor} \frac{1}{j}$ which he introduced in this context. A generalization of his result is proved in the case of composite $N$, and a new proof given of a sharpened result by Skula (2008) \cite{Skula2008} when $N$ is even. We also sharpen the criteria given by Emma Lehmer in 1938 for a Wieferich prime to be simultaneously a Mirimanoff prime.

\noindent
\textit{Keywords}: Fermat quotient, Wieferich prime, Mirimanoff prime
\end{abstract}

\section{Introduction}

\fbox {
\parbox{\linewidth}{
Note: In the version of this paper uploaded to the Mathematics arXiv on 6 July 2014, I noted that it had just been brought to my attention by a correspondent that an article by Dilcher \& Skula published in 2011 \cite{DilcherSkula2011} contained many of the same findings. However, I do not believe that a retraction of the present piece is in order. With no intended disrespect to those authors, I must point out that a recent consultation of my files reminded me that while this paper was only first uploaded to the Mathematics arXiv in March 2011, it had been available between 20 June 2009 and 31 March 2011 under the same title and in substantially the same form on my personal website \cite{Dobson2009}. While the link to the file was removed from my website when the paper was first uploaded to the Mathematics arXiv, and the version on my website was apparently never captured by the Internet Archive, the file itself (which was not subsequently altered) remains on the University of Winnipeg Library's web server, and the system datestamp corroborates the upload date.
}
}

\bigskip

\noindent
For the Fermat quotient $q_p(b) = (b^{p - 1} - 1)/p \pmod{p}$ we employ wherever possible one of the briefer notations $q(b)$ or $q_b$. Henceforth all congruences are assumed to be mod $p$ unless otherwise stated, and $\lfloor \cdot \rfloor$ signifies the greatest-integer function.

The fact that Fermat quotients can be expressed as sums involving reciprocals of integers in $\left\{ 1, p-1\right\}$ was discovered in 1850 for the case $b = 2$ by Eisenstein, who gives

\begin{displaymath}
2 \cdot q_2 \equiv 1 - \frac{1}{2} + \frac{1}{3} - \frac{1}{4} + \ldots + \frac{1}{(p - 2)} - \frac{1}{(p - 1)}.
\end{displaymath}

\noindent
Subsequent researches in this direction have focused on developing equivalent results entailing fewer terms. Later Wieferich and Mirimanoff demonstrated the connection of Fermat quotients with the first case of Fermat's Last Theorem (FLT). This connection retains its historical interest despite the full proof of FLT by Wiles in 1995, and Skula's demonstration in 1992 \cite{Skula1992} that the failure of the first case of FLT would imply the vanishing of many similar sums but with much smaller ranges, which cannot be evaluated in terms of Fermat quotients.

Lerch's remarkable formula of 1905 (\cite{Lerch}, p.\ 474) is

\begin{equation} \label{eq:Lerch}
N \cdot q_p(N) \equiv \sum_{k=1}^{N-1} k \cdot s(k, N),
\end{equation}

\noindent
where

\begin{displaymath}
s(k, N) = \sum_{\substack{j=\lfloor\frac{kp}{N}\rfloor + 1\\ j \neq p}}^{\lfloor\frac{(k + 1)p}{N}\rfloor}
\frac{1}{j},
\end{displaymath}

\noindent
and it is always assumed that $p$ is sufficiently large that $s(k, N)$ contains at least one element; the provision $j \neq p$ is necessary when $k+1=N$. For the background of this result see the historical study by Lepka \cite{Lepka2000}, and for a detailed exposition in English, Agoh \cite{Agoh1996}, and Agoh et al.\ (\cite{AgohDilcherSkula}, pp.\ 32--35). In the case of composite $N$, the left-hand side of (\ref{eq:Lerch}) is evaluated using Eisenstein's logarithmetic property $q(ab) = q(a) + q(b)$. The fact that the terms in $\left\{\frac{(p+1)}{2}, p-1\right\}$ (descending) are the complements $\pmod{p}$ of those in $\left\{ 1, \frac{(p-1)}{2}\right\}$ (ascending) implies that

\begin{equation} \label{eq:Complementarity}
s(k, N) \equiv -s(N-1-k, N).
\end{equation}

\noindent
For ease of comparison with previous literature, in our final results we usually restrict $k$ so as to be less than $\frac{(p-1)}{2}$, or to be of a particular parity; but in the proofs we use whichever form seems more intelligible or expressive in the given situation. 

As Lerch himself noted, the complementarity of the terms about the middle of the range $\left\{ 1, p-1\right\}$ leads to considerable simplification of (\ref{eq:Lerch}) above, with $s\left(\frac{(N-1)}{2}, N\right)$ vanishing for odd $N$, so that

\begin{subequations}

\begin{multline} \label{eq:LerchSplit}
N \cdot q_p(N) \equiv \\
-(N-1) \cdot s(0, N) - (N-3) \cdot s(1, N) - (N-5) \cdot s(2, N) - \ldots \\
- s\left(\frac{N}{2}-1, N\right) \quad [N \textrm{ even}] \quad
\end{multline}

\begin{multline} 
N\cdot q_p(N) \equiv \\
-(N-1)\cdot s(0, N) - (N-3) \cdot s(1, N) - (N-5) \cdot s(2, N) - \ldots \\
- 2 \cdot s\left(\frac{(N-1)}{2} - 1, N\right) \quad [N \textrm{ odd}]. \quad
\end{multline}

\end{subequations}

Lerch's results for the Fermat quotient are so comprehensive as to subsume all those previously achieved. The cases $N=2$ and $N=4$ give, respectively, Glaisher's 1901 results (\cite{Glaisher}, pp.\ 21-22, 23) $s(0, 2), \equiv -2 \cdot q_2$ and $s(0, 4) \equiv -3 \cdot q_2$. The case $N=3$ gives another result of Glaisher, that $s(0, 3) \equiv -\frac{3}{2} \cdot q_3$ (\cite{Glaisher3}, p.\ 50). The unsatisfactory account of the literature on the Fermat quotient given in Dickson's \textit{History} \cite{Dickson} credits the results for $N=2$ and $N=4$ to later writers (1:111, n. 32a; 1:109, n. 21), and incorrectly implies (1:109, n. 18) that Glaisher only proved the case of $N=3$ for $p \equiv 1 \pmod 3$, when in fact he also proved it for $p \equiv 2 \pmod 3$. Lerch himself (\cite{Lerch}, p.\ 476, equations 14 and 15) pretty much explicitly writes out $2 \cdot s(0, 5) + s(1, 5) \equiv -\frac{5}{2} \cdot q_5$, correcting work of Sylvester (\cite{Sylvester}, pp.\ 161--162). Lerch could have used his formulae to evaluate $s(0, 6)$, which would have completed the treatment of the cases where $k|6$, and so it is convenient to group these instances among the ``classical'' results which were completely settled and systematized by Lerch's method. However, in actuality the famous 1938 paper by Emma Lehmer (\cite{Lehmer1938}, pp.\ 356ff), which points out that $s(0, 6) \equiv -2 \cdot q_2 - \frac{3}{2} \cdot q_3$, seems to be the first source ever to give the evaluation explicitly. Likewise, we have not seen in any publication prior to that of Granville and Sun in 1996 (\cite{GranvilleSun}, p.\ 136) an explicit statement that the two instances of $s(k, 12)$ lying at the center of the range $\left\{ 1, \frac{(p - 1)}{2}\right\}$ can be evaluated by subtraction from known results, giving $s(2, 12) \equiv -q_2 + \frac{3}{2} \cdot q_3$, and $s(3, 12) \equiv 3 \cdot q_2 - \frac{3}{2} \cdot q_3$. However, as it happens, this evaluation is supplied by Vandiver's corollary of 1917 (see our equation \ref{eq:VandiverDifference} below), and so also fits naturally among the ``classical'' results, which are summarized in Table 1 below.

Some of these results have been extended to higher moduli. Lehmer (\cite{Lehmer1938}, p. 305, eq.\ 45) gives $s(0, 2) \equiv q_2 \equiv -q_2 + pq_2^2$ (cf. \cite{Cai}). Zhi-Hong Sun (\cite{ZHSun2000}, p. 208) shows $s(0, 2) \equiv -2q_2 + pq_2^2 -\frac{2}{3}p^2q_2^3 - \frac{7}{12}p^2B_{p-3} \pmod{p^3}$. Elsewhere the same author (\cite{ZHSun2008}, pp.\ 289, 298) completely evaluates every instance of $s(k, 4) \pmod{p^3}$, and supplies sketches of similar evaluations for $s(k, 3)$ and $s(k, 6)$, which are written out explicitly in (\cite{ZHSun2012}, p.\ 215). These evaluations are reproduced in Tables 2 and 3 below.

Dilcher \& Skula (\cite{DilcherSkula1995}, p.\ 389) report numerical investigations of all possible instances of $s(k, N)$ with $N \le 46$ for the two values of $p$ for which $q_2$ vanishes, namely 1093 and 3511. Because another portion of their paper contains computational errors (reported below), we repeated these calculations, and found them to be correct; furthermore, we examined all possible cases of $N$, such that $s(k, N)$ contained at least three terms. (As Dilcher \& Skula note, other than the central terms in the cases $s\left(\frac{(N-1)}{2}, N\right)$ with odd $N$ which trivially vanish and therefore require no testing, no $s(k, N)$ can contain two adjacent terms that are complementary mod $p$, but three is a sharp minimum as $s(0, 3)$ contains only three terms in the case $p$ = 11, where it vanishes.) The vanishing of $q_2$ indeed occurs only for the ``classical'' cases, eliminating the possibility that any other sums involving consecutive terms could be simple multiples of $q_2$. We therefore conclude that our Theorems 1 and 2 comprise essentially all the linear relations which pertain among sums of Lerch's type.

\section{Supplementary Notations}

Certainly, not all sums figuring in the literature of Fermat quotients can be reduced to Lerch's type, least of all those containing the numbers of Bernoulli, Euler, Fibonacci, Lucas, or Pell. However, the sums studied here are the simplest representatives of an important family of interrelated sums whose other members we designate as follows: 

\begin{multline} \label{eq:sPrime}
s^\prime(k, N) = \textrm{terms with odd denominators in } s(k, N)\\
\equiv -\frac{1}{2} \cdot s(N-1-k, 2N) \equiv \frac{1}{2} \cdot s(N+k, 2N) 
\end{multline}

\begin{multline} \label{eq:sDoublePrime}
s{''}(k, N) = \textrm{terms with even denominators in } s(k, N) \equiv \frac{1}{2} \cdot s(k, 2N)\\
\end{multline}

\begin{multline} \label{eq:sTriplePrime}
s{'''}(k, N) = \textrm{terms with denominators divisible by 3 in } s(k, N) \equiv \frac{1}{3} \cdot s(k, 3N)\\
\end{multline}

\begin{multline} \label{eq:sStar}
s^{\ast}(k, N) = s{''}(k, N) - s^\prime(k, N) \equiv \frac{1}{2} \cdot s(k, 2N) - \frac{1}{2} \cdot s(N + k, 2N)\\
\end{multline}

\begin{multline} \label{eq:K}
K(r, N) = \textrm{terms in } s(0, 1) \textrm{ with denominators congruent to } rp \bmod N\\
\equiv \frac{1}{N} \cdot s(N-r, N) \equiv -\frac{1}{N} \cdot s(r-1, N)
\end{multline}

The evaluation of $s^\prime(k, N)$ (equation \ref{eq:sPrime}) is a simple consequence of the fact that a series expressed as a sum of terms of odd denominator may be condensed into a smaller range of unrestricted terms, as follows: 

\begin{equation} \label{eq:condensation}
\sum_{\substack{1 \\ j \textrm{ odd}}}^{n} \frac{1}{j} \equiv
\sum_{\substack{1 \\ j \textrm{ odd}}}^{2\lfloor(n-1)/2\rfloor + 1} \frac{1}{j} \equiv
- \sum_{\substack{2\lfloor(p-n+1)/2\rfloor \\ j \textrm{ odd}}}^{p-1} \frac{1}{j} \equiv
-\frac{1}{2} \sum_{\substack{\lfloor(p-n+1)/2\rfloor \\}}^{(p-1)/2} \frac{1}{j}.
\end{equation}

\noindent
Conversely, a series expressed as a sum only of terms of even denominator may be simplified merely by factoring out the 2 in the denominator of its summand, whence the formula for $s{''}(k, N)$ (equation \ref{eq:sDoublePrime}). Stern (\cite{Stern}, p.\ 184) gives $s^\prime(0, 1) \equiv q_2$, and Zhi-Hong Sun (\cite{ZHSun2000}, p.\ 210, Remark 5.3) extends this result to $p^3$. Reformulating work of Glaisher (\cite{Glaisher}, p. 6), Vandiver (\cite{Vandiver}, pp.\ 111--12) shows that $s^\prime(0, 1) + \frac{1}{2}p\cdot \left\{ s^\prime(0, 1) \right\}^2 \equiv q_2 \pmod{p^2}$. The formula for $s^{\ast}(k, N)$ (\ref{eq:sStar}) follows immediately from (\ref{eq:sPrime}) and (\ref{eq:sDoublePrime}), as does the fact that

\begin{equation}
s^{\ast}(k, N) \equiv s^{\ast}(N-1-k, N).
\end{equation}

\noindent
It is easy to prove the following generalization of a statement in Sun and Sun (\cite{Sun+Sun1992}, p.\ 385):

\begin{equation} \label{eq:sStar2}
s(k, N) + s^{\ast}(k, N) \equiv s(k, 2N),
\end{equation}

\noindent
and an alternate form thereof:

\begin{equation} \label{eq:sStar3}
s^{\ast}(k, N) - s(k, N) \equiv s(N-1-k, 2N).
\end{equation}

\noindent
An important special case of this function, derivable from Corollary 1 below (\ref{eq:Corollary1}), is

\begin{equation} \label{eq:0_N}
s^{\ast}(0, N) \equiv -s(1, 2N),
\end{equation}

\noindent
and this in turn gives

\begin{equation} \label{eq:1_N}
s^{\ast}(1, N) \equiv s(1, 2N) - s(1, N).
\end{equation}

\noindent
It is clear from the definition that for odd $N$,

\begin{equation} \label{eq:StarOdd}
s^{\ast}\left(\frac{N-1}{2}, N\right) \equiv s\left(\frac{N-1}{2}, 2N\right),
\end{equation}

\noindent
while for $N$ oddly even,

\begin{equation} \label{eq:StarOddlyEven}
s^{\ast}\left(\frac{N-2}{2}, N\right) \equiv \frac{1}{2} \cdot s\left(\frac{N-2}{4}, N\right).
\end{equation}

\noindent
Some specific values of this function are well known: $s^{\ast}(0, 1) \equiv -s(1, 2) \equiv -2 \cdot q_2$ is the classic result from 1850 of Eisenstein (\cite{Eisenstein}, p.\ 41), and $s^{\ast}(0, 2) \equiv -s(1, 4) \equiv -q_2$ is given by Stern (\cite{Stern}, pp.\ 185--86). Glaisher's statement of the latter result (\cite{Glaisher}, pp.\ 22--23) is a rediscovery and does not warrant the level of credit assigned in Dickson's \textit{History} \cite{Dickson}, 1:111. Zhi-Hong Sun (\cite{ZWSun2011}, p.\ 2515) gives $s^{\ast}(0, 2) \equiv -q_2 -\frac{1}{2}p \cdot q_2^2 + (-1)^{(p+1)/2}pE_{p-3} \pmod{p^2}$, where $E$ is an Euler number. $K(r, N)$ (equation \ref{eq:K}), the proof of the formula for which is defered to (\ref{eq:collectArrayOne}) below, corresponds to the $K_m(s, p)$ of Zhi-Hong Sun \cite{ZHSun1992} and to the $K_p(r, m)$ of Zhi-Wei Sun \cite{ZWSun2002}; we however omit the parameter $p$ to simplify the notation and make it more uniform with the rest. Obviously $K(0, 2) \equiv s{''}(0, 1)$ and $K(1, 2) \equiv s^\prime(0, 1)$. Zhi-Hong Sun (\cite{ZHSun2008}, pp.\ 281, 286--288, 303) evaluates ${K}(1, 3)$, ${K}(1, 4)$, $K(1, 6)$, and $K(-1, 4)$ mod $p^3$, and $K(-1, 3)$ mod $p^2$, as well as some analogous sums over smaller ranges.

The function $s{'''}(k, N)$ (\ref{eq:sTriplePrime}) is not central to our argument, and for now we state without proof the only results for it that can be expressed solely in terms of Fermat quotients:

\begin{align}
s{'''}(0, 1) & \equiv \frac{1}{3} \cdot s(0, 3)  \;\: \equiv             -\frac{1}{2} \cdot q_3 \\
s{'''}(0, 2) & \equiv \frac{1}{3} \cdot s(0, 6)  \;\: \equiv             -\frac{2}{3} \cdot q_2 - \frac{1}{2} \cdot q_3 \\
s{'''}(1, 2) & \equiv \frac{1}{3} \cdot s(1, 6)  \;\: \equiv \hphantom{-} \frac{2}{3} \cdot q_2 \\
s{'''}(2, 4) & \equiv \frac{1}{3} \cdot s(2, 12)      \equiv             -\frac{1}{3} \cdot q_2 + \frac{1}{2} \cdot q_3 \\
s{'''}(3, 4) & \equiv \frac{1}{3} \cdot s(3, 12)      \equiv \hphantom{-} \quad\;\:           q_2 - \frac{1}{2} \cdot q_3.
\end{align}

\noindent
All of these have obvious relevance to Lehmer's problem (see section \ref{Lehmersproblem} below). For some variations on $s{'''}(0, 2)$ and $s{'''}(1, 2)$ mod $p^3$ see Zhi-Hong Sun (\cite{ZHSun2008}, p.\ 288). 

It should be noted that Glaisher \cite{Glaisher4} gave evaluations of sums of reciprocals of integers restricted to each residue class for the moduli 2, 3, 4, and 6, in terms of Bernoulli numbers.

\section{Alternate expressions for Lerch's sums}

\subsection{Harmonic numbers} \label{Harmonic_numbers}

\noindent
It will be apparent from the definitions above that

\begin{equation} \label{eq:HarmonicSeries}
\sum_{j=0}^{k} s(j, N) = \sum_{j=1}^{\lfloor\frac{(k + 1)p}{N}\rfloor} \frac{1}{j},
\end{equation}

\noindent
and that in the simplest case, where $k=0$, the right-hand side is nothing but a partial sum of the Harmonic series, $H_g = 1 + \frac{1}{2} + \frac{1}{3} + \dots + \frac{1}{g}$, where $g = \lfloor\frac{p}{N}\rfloor$. In this notation, Lerch's congruence (\ref{eq:Lerch}) assumes the form in which it is given by Vandiver (\cite{Vandiver}, p.\ 114), which will be considered in greater detail in section \ref{Vandiver_corollary} below:

\begin{equation} \label{eq:LerchHarmonic}
\sum_{j=1}^{N-1} H_{\lfloor jp/N \rfloor} \equiv -N \cdot q_N.
\end{equation}

\noindent
Letting $N = p$ gives $\sum_{j=1}^{p-1} H_j \equiv 1 \pmod{p}$, though this result is weak compared with the congruence $\sum_{j=1}^{p-1} H_j \equiv 1 - p \pmod{p^3}$ given by Z.\,W. Sun (\cite{ZWSun2011b}, p.\ 418). But as a generalization of (\ref{eq:Complementarity}) we have $H_{p-k} \equiv H_{k-1} \pmod{p}$, so the equivalents in Harmonic numbers of Lerch's short sums (\ref{eq:LerchSplit}) are

\begin{subequations}

\begin{equation} \label{eq:LerchHarmonicSplit}
\sum_{j=1}^{N/2} H_{\lfloor jp/N \rfloor} \equiv -\frac{1}{2} N \cdot q_N - q_2 \quad [N \textrm{ even}]
\end{equation}

\begin{equation} 
\sum_{j=1}^{(N-1)/2} H_{\lfloor jp/N \rfloor} \equiv -\frac{1}{2} N \cdot q_N \quad [N \textrm{ odd}; N \neq p].
\end{equation}

\end{subequations}

\noindent
In the first row, letting $N = p-1$ gives $\sum_{j=1}^{(p-1)/2} H_j \equiv \frac{1}{2} - q_2 \pmod{p}$. Further congruences in this vein will be found in \cite{ZWSun2011b} and \cite{Mestrovic2014}.

Clearly, in view of the results of Glaisher previously cited, and as noted in \cite{Agoh1992}, pp.\ 90--91, we have $H_{\lfloor p/4 \rfloor} \equiv -3 \cdot q_2$, $H_{\lfloor p/3 \rfloor} \equiv -\frac{3}{2} \cdot q_3$, $H_{\lfloor p/2 \rfloor} \equiv -2q_2 \pmod{p}$; also, $H_{p-1} \equiv 0 \pmod{p^2}$ $(p>3)$ by Wolstenholme's theorem. Wieferich primes are thus characterized by the condition $H_{(p-1)/2} \equiv H_{\lfloor p/4 \rfloor} \equiv 0 \pmod{p}$, Mirimanoff primes by the condition $H_{\lfloor 2p/3 \rfloor} \equiv H_{\lfloor p/3 \rfloor} \equiv 0 \pmod{p}$, and both by the apparently rare condition $H_{\lfloor 2p/N \rfloor} \equiv H_{\lfloor p/N \rfloor} \equiv 0 \pmod{p}$ for some $N>1$. Primes for which $p$ divides only $H_{p-1}$, $H_{p^2-p}$, and $H_{p^2-1}$ are known as the harmonic primes, and the exceptions as the anharmonic primes (\cite{EswarathasanLevine}, \cite{Boyd}). It should be noted that $\mathcal{H}_{p, m}(p-1)$ in the important paper by Hao Pan \cite{Pan2011} agrees with our $K(1, m)$, but only to the modulus $p$. We cannot claim to have made a thorough survey of the immense literature on Harmonic sums, and it is therefore possible that some results pertinent to the present paper may have been overlooked.

The partial odd Harmonic series $H_g^\prime = 1 + \frac{1}{3} + \frac{1}{5} + \dots + \frac{1}{g}$, where $g$ is the greatest odd integer $\leq \frac{p}{N}$, also has an historical connection with the study of the Fermat quotient; and $H_{p-1}^\prime$ is the same as the sum designated $s^\prime(0, 1)$ elsewhere in the present paper. As previously noted, Stern \cite{Stern} gives $H_{p-1}^\prime \equiv q_2 \pmod{p}$, and Zhi-Hong Sun \cite{ZHSun2000} has extended this result to $p^3$. Glaisher (\cite{Glaisher}, pp.\ 22, 26) gives $H_{p-1}^\prime \equiv s^{\ast}(0, 2) \equiv -\frac{1}{2} s(0, 2) \pmod{p}$ and

\begin{displaymath}
q_2 \equiv H_{p-1}^\prime + \frac{1}{2}p\cdot (H_{p-1}^\prime)^2 \pmod{p^2},
\end{displaymath}

\noindent
so that $q_2 \equiv 0 \pmod{p^2}$ if and only if $H_{p-1}^\prime \equiv 0 \pmod{p^2}$; a much more succinct proof of this result appears in Vandiver (\cite{Vandiver}, pp.\ 111--12).

\subsection{Binomial coefficients} \label{binomial_coefficients}

\noindent
Another equivalence given in Lehmer's remarkably fertile paper (\cite{Lehmer1938}, p.\ 360) is

\begin{displaymath}
H_k \equiv \frac{(-1)^{k+1} \binom{p-1}{k} + 1}{p} \pmod{p},
\end{displaymath}

\noindent
whence

\begin{displaymath}
s(0, N) := H_{\lfloor p/N \rfloor} \equiv \frac{(-1)^{\lfloor p/N \rfloor+1} \binom{p-1}{\lfloor p/N \rfloor} + 1}{p} \pmod{p}.
\end{displaymath}

\noindent
Because $\binom{p-1}{\lfloor kp/N \rfloor} \equiv (-1)^{\lfloor kp/N \rfloor + 1}\left(pH_{\lfloor kp/N \rfloor} - 1 \right)$ $\pmod{p^2}$, products taken over binomial coefficients of this kind simplify mod $p$ to sums over $H_{\lfloor kp/N \rfloor}$. Thus it is straighforward to verify that Lerch's main theorem as applied to Harmonic numbers (\ref{eq:LerchHarmonic}) is a consequence of Granville's celebrated theorem on products of binomial coefficients (\cite{Granville}, pp. 257--58):

\begin{displaymath}
(-1)^{\frac{p-1}{2}(N-1)} \prod_{j=1}^{N-1} \binom{p-1}{\lfloor jp/N \rfloor} \equiv N^p - N + 1 \pmod{p^2}.
\end{displaymath}

\noindent
Zhi-Wei Sun has successively sharpened Granville's result, first showing \cite{ZWSun2001} that

\begin{equation} \label{eq:SunSplit}
\begin{split}
& (-1)^{\frac{p-1}{2}\lfloor N/2 \rfloor} \prod_{j=1}^{\lfloor N/2 \rfloor} \binom{p-1}{\lfloor jp/N \rfloor} \\
& \equiv \begin{dcases}
\left(\frac{N}{p}\right)\left(\frac{2}{p}\right) \left\{1 + \frac{N}{2}\left(N^{p-1} - 1\right) + \left(2^{p-1} - 1\right) \right\} \pmod{p^2} & \text{if $2 \mid N$} \\
\left(\frac{N}{p}\right) \left\{1 + \frac{N}{2}\left(N^{p-1} - 1\right) \right\} \pmod{p^2} & \text{if $2 \nmid N$}
\end{dcases}
\end{split}
\end{equation}

\noindent
where we have simplified the notation. It can be verified without difficulty that these give (\ref{eq:LerchHarmonicSplit} \& b). Later \cite{ZWSun2006}, in the cases where $N$ is not a perfect square, he was able to restrict the product such that $j$ runs over only the quadratic residues, or only the non-residues, of $N$. Since the counts of the residues or nonresidues is always $\le N/2$, these congruences refine (\ref{eq:LerchHarmonicSplit} \& b). The pertinent congruence is in fact sharp enough to resolve the value of $s(0,5)$; this case is worked out in the review of this paper in \textit{Zentralblatt MATH Database}, though with a typographical error in the index of the Fibonacci number, which should not contain a factor of 3. The case of $s(k, 24)$ is considered in section \ref{cases_24} below.

The study of binomial coefficients of the type discussed in this section is enormously extended in \cite{KuzumakiUrbanowicz}, where the treatment of the case $N$ = 24 enables the complete evaluation of $s(k, 24)$ (see section \ref{cases_24} below).

\subsection{Bernoulli numbers}

Though (\ref{eq:HarmonicSeries}) could be written out in full generality in terms of Bernoulli polynomials and Bernoulli numbers, it is only the case of $k=0$ that seems worth mentioning here. Then, an important theorem of Emma Lehmer (\cite{Lehmer1938}, p.\ 351, eq.\ 8, setting $2k-1 = p-2$), gives

\begin{equation} \label{eq:HarmonicSeries2}
\begin{split}
s(0, N) := \sum_{j=1}^{\lfloor\frac{p}{N}\rfloor} \frac{1}{j} & \equiv \frac{1}{p-1} \left\{ B_{p-1}\left(\frac{s}{n}\right) - B_{p-1} \right\} - \frac{p}{n} B_{p-2}\left(\frac{s}{n}\right) \pmod{p^2} \\
& \equiv -B_{p-1}\left(\frac{s}{n}\right) + B_{p-1} \pmod{p},
\end{split}
\end{equation}

\noindent
where $s$ is the least positive residue of $p \pmod{N}$. Lehmer's discussions of the cases $N = 2, 3, 4, 6$, and that by Granville \& Sun \cite{GranvilleSun} of the cases $N = 5, 8, 10, 12$, are not in most instances the first in the literature, but these writings may be particularly recommended for the coherence of their methods.

\section{A partial generalization of Lerch's formula for composite $N$}

If $N$ is prime, then Lerch's formula (\ref{eq:Lerch}) appears to be the only linear relation which holds among his sums of order $N$. However, if $N$ is composite, (\ref{eq:Lerch}) can be manipulated in a way which does not appear to have been previously stated. For example, if $N$ is even, we have 

\begin{equation} \label{eq:LerchDiff}
\begin{split}
N \cdot q_p(N) & - 2N \cdot q_p\left(\frac{N}{2}\right) \equiv N \cdot q_p(2) \\
& \equiv -\left\{ s(0, N) + s(2, N) + s(4, N) + \ldots + s(N - 2, N)\right\}.
\end{split}
\end{equation}

\noindent
This coincides with a result given below, but because we want to prove a generalization of the formula (when $N$ is composite) we cannot work directly from (\ref{eq:Lerch}), but rather must develop a somewhat broader version of the underlying theory. 

The essential idea of Lerch's formula is that for fixed $N$, every pair of values of $s(k, N)$ is connected with every other value through the sharing of multiples of each other's terms $\pmod{p}$. Specifically, $s(k, N)$ by definition contains the reciprocals of all $j$ in the range $\left\{\lfloor\frac{(k + 1)p}{N}\rfloor, \lfloor\frac{kp}{N}\rfloor + 1\right\}$, and if $M$ is any divisor of $N$ (including $N$ itself), and $r$ the residue of $k \pmod{M}$, then the corresponding values of $Mj$ reduced $\bmod p$ are distributed among the sums $s(rM, N)$ through $s(rM+M-1, N)$ as follows: 

\begin{tabular} { l r }
\textbf{range of $k$}                                       & \textbf{residue of $Mj \bmod{p}$} \\
\hline
$\left\{ 0, \frac{N}{M} - 1 \right\}$                       & $0$ \\
$\left\{ \frac{N}{M}, 2\cdot\frac{N}{M} - 1 \right\}$       & $-p$ \\
$\left\{ 2\cdot\frac{N}{M}, 3\cdot\frac{N}{M} - 1 \right\}$ & $-2p$ \\
$\ldots$                                                    & $\ldots$ \\
$\left\{ \left( M-1 \right) \frac{N}{M}, N-1 \right\}$      & $-(M-1)p$ \\
\end{tabular}

\noindent
That is, the values of $Mj$ in the first row fall in $\left\{ 1, p-1\right\}$, while those in the second row fall in $\left\{ p+1, 2p-1\right\}$ and must be reduced by $p$, those in the third row fall in $\left\{ 2p + 1, 3p - 1\right\}$ and must be reduced by $2p$, etc. 

We can evaluate the sum of this two-dimensional array of terms in two different ways. First, we can collect the terms belonging to each value of $Mj \bmod{p}$ as follows:

\begin{multline*}
s(0, M)  \equiv M \times \left\{ \textrm{terms } \equiv 0 \pmod{M} \textrm{ in } s(0, 1)\right\}  \equiv M \cdot K(0, M) \qquad \\ 
\end{multline*}

\begin{multline*}
s(1, M)  \equiv M \times \left\{ \textrm{terms } \equiv -p \pmod{M} \textrm{ in } s(0, 1)\right\}  \equiv M \cdot K(-1, M) \qquad \\ 
\end{multline*}

\begin{multline*}
s(2, M)  \equiv M \times \left\{ \textrm{terms } \equiv -2p \pmod{M} \textrm{ in } s(0, 1)\right\}  \equiv M \cdot K(-2, M) \qquad \\
\end{multline*}

\begin{multline*}
\ldots \ldots \dots
\end{multline*}

\begin{multline} \label{eq:collectArrayOne}
s(M-1, M) \equiv M \times \left\{ \textrm{terms } \equiv -(M-1)p \pmod{M} \textrm{ in } s(0, 1)\right\} \\
\equiv M \cdot K(-(M-1), M).
\end{multline}

\noindent
In other words $N \cdot K(r, N) \equiv s(-r, N) \equiv -s(r-1, N)$, a result which coincides with one of Zhi-Hong Sun (\cite{ZHSun1992}, pt. 3, p.\ 90, Corollary 3.1). Secondly, we can collect the terms belonging to each value of $r$, as follows:

\subsubsection*{Theorem 1}

\begin{multline*}
s(0, N) + s\left(\frac{N}{M}, N\right) + s(2\cdot\frac{N}{M}, N) + \ldots + s\left(\left(M - 1\right)\frac{N}{M}, N\right) \\
\equiv M\left\{ s(0, N) + s(1, N) + s(2, N) + \ldots + s(M - 1, N)\right\} \\
\equiv M \cdot s(0, N/M) \qquad
\end{multline*}

\begin{multline*}
s(1, N) + s\left(\frac{N}{M} + 1, N\right) + s(2\cdot\frac{N}{M}+1, N) + \ldots + s\left(\left(M-1\right)\frac{N}{M}+1, N\right) \\
\equiv M\left\{ s(M, N) + s(M + 1, N) + s(M + 2, N) + \ldots + s(2M - 1, N)\right\} \\
\equiv M \cdot s(1, N/M) \qquad
\end{multline*}

\begin{multline*}
\ldots \ldots \ldots
\end{multline*}

\begin{multline} \label{eq:collectArrayTwo}
s\left(\frac{N}{M}-1, N\right) + s(2\cdot\frac{N}{M} - 1, N) + s(3\cdot\frac{N}{M}-1, N) + \ldots + s(N-1, N) \\
\equiv M\left\{ s(N-M, N) + s(N-M+1, N) + s(N-M+2, N) + \ldots + s(N-1, N)\right\} \\
\equiv M \cdot s\left(\frac{N}{M}-1, N/M\right).
\end{multline}

\noindent
Now, while all the relations in (\ref{eq:collectArrayTwo}) are valid, those in the bottom half of the display are by symmetry simply the negatives of those in the opposite rows in the top half, while if $M$ is odd the statement made in the middle row is truistic because each side consists of sums of terms symmetrically distributed about the ``center line'' and thus vanishes mod $p$. Furthermore, these relations are degenerate when $M = N$, so $M$ must be a proper divisor of $N$, and thus $N$ cannot be prime. 

The relationships expressed in (\ref{eq:collectArrayTwo}) may be employed in two different ways. The simpler expressions at the far right, if tractable, will provide an explicit evaluation of the sums on the far left, although as we shall see in the cases of $N = 9$ and $N = 18$, it is not a foregone conclusion that the sums are better known for smaller values of $N$. Alternatively, the sums on the far left and in the middle are homogeneous, and the congruence always admits of some simplification, both by direct cancelation and by use of the rule $s(k, N) \equiv -s(N-k-1, N)$. 

We now examine some of the implications of (\ref{eq:collectArrayTwo}) for particular cases of $M$, reserving a discussion of particular cases of $N$ until the end.

\subsection{The case $M=2$}

The case $M=2$ of (\ref{eq:collectArrayTwo}) is especially interesting. Let $x = N/M$; then:

\begin{equation} \label{eq:M=2}
\begin{tabular} {@{} l @{ } l @{ } l @{}}
$s(0, 2x) + s(x, 2x)$      & $\equiv 2\left\{ s(0, 2x) + s(1, 2x)\right\}$      & $\equiv 2 \cdot s(0, x)$ \\
$s(1, 2x) + s(x+1, 2x)$    & $\equiv 2\left\{ s(2, 2x) + s(3, 2x)\right\}$       & $\equiv 2 \cdot s(1, x)$ \\
\ldots                     & \ldots                                             & \ldots \\
$s(x-1, 2x) + s(2x-1, 2x)$ & $\equiv 2\left\{ s(2x-2, 2x) + s(2x-1, 2x)\right\}$ & $\equiv 2 \cdot s(x-1, x)$, \\
\end{tabular}
\end{equation}

\noindent
with the usual redundancy in the second half of the display. The first row of (\ref{eq:M=2}) yields an important corollary which will be frequently invoked in the following pages:

\subsubsection*{Corollary 1}

\begin{multline} \label{eq:Corollary1}
s(0, 2x) + s(x, 2x) \equiv 2\left\{ s(0, 2x) + s(1, 2x)\right\} \\
\Rightarrow s(0, 2x) + 2 \cdot s(1, 2x) - s(x, 2x) \equiv 0 \\
\Rightarrow s(0, 2x) + 2 \cdot s(1, 2x) + s(x-1, 2x) \equiv 0 
\end{multline}

\noindent
Among relations connecting the sums $s(k, N)$ for even $N$ which do not depend on knowledge of the value of some sum with smaller $N$, this one is noteworthy in that it relates only three terms, whereas Lerch's (\ref{eq:LerchSplit}) involves $\lfloor N/2\rfloor$ terms. Other relations compare sums of two values of $s(k, N)$ for even $N$ with one or two values of $s(k, N/2)$. In (\ref{eq:Corollary1}), subtract both sides from $2s(0, x) + s\left(\frac{x}{2}-1, x\right)$. Then

\begin{equation}
s(0, 2x) + s(x-2, 2x) \equiv 2s(0,x) + s\left(\frac{x}{2}-1, x\right). 
\end{equation}

\noindent
In (\ref{eq:Corollary1}), subtract $s(0, x)$ from both sides. Then

\begin{equation}
s(1, 2x) + s(x-1, 2x) \equiv -s(0,x). 
\end{equation}

\noindent
Also, using (\ref{eq:Corollary1}) and (\ref{eq:sStar}), 

\begin{equation} \label{eq:sStar_rule}
s^{\ast}(0, x) \equiv -s(1, 2x). 
\end{equation}

\subsubsection*{Corollary 2}

Adding together the first two rows of (\ref{eq:M=2}), with the use of (\ref{eq:sStar}) and some manipulation we derive the relation: 

\begin{equation}
s^{\ast}(0, x) + s^{\ast}(1, x) \equiv -s(1, x). 
\end{equation}

\noindent
This relation, which corresponds to a well-known Harmonic-number identity (see for example \cite{Larsen}, p.\ 187, Remark 15.21), accounts for the equivalence of some evaluations given by Zhi-Hong Sun (\cite{ZHSun1992}, pt. 3, Theorem 3.2, nos.\ 2 and 3).

\subsection{The case $M=N/2$}

In the first two rows of (\ref{eq:collectArrayTwo}), let $N$ be even, and $M=N/2$. Then 

\subsubsection*{Corollary 3}

\begin{subequations}

\begin{multline} \label{eq:Corollary3}
s(0, N) + s(2, N) + s(4, N) + \ldots + s(N - 2, N) \\
\equiv \frac{N}{2}\left\{ s(0, N) + s(1, N) + s(2, N) + \ldots + s\left(\frac{N}{2} - 1, N\right)\right\} 
\\ \equiv \frac{N}{2} \cdot s(0, 2)  \equiv -N \cdot q_2 \qquad
\end{multline}

\begin{multline} 
s(1, N) + s(3, N) + s(5, N) + \ldots + s(N - 1, N) \\
\equiv \frac{N}{2}\left\{ s\left(\frac{N}{2}, N\right) + s\left(\frac{N}{2} + 1, N\right) + s\left(\frac{N}{2} + 2, N) + \ldots + s(N - 1, N\right)\right\} \\
\equiv \frac{N}{2} \cdot s(1, 2)  \equiv N \cdot q_2.
\end{multline}

\end{subequations}

\noindent
Now let $N=p-1$; then the sums $s(k, p-1)$ divide $\left\{ 1, p-1\right\}$ into equal pieces of length 1, so that:

\begin{subequations}
\begin{align}
\frac{1}{1} + \frac{1}{3} + \frac{1}{5} + \ldots + \frac{1}{(p - 2)} & \equiv s^\prime(0, 1)  \equiv -\frac{1}{2} \cdot s(0, 2) \equiv q_2 \\
\frac{1}{2} + \frac{1}{4} + \frac{1}{6} + \ldots + \frac{1}{(p - 1)} & \equiv s{''}(0, 1)  \equiv -\frac{1}{2} \cdot s(1, 2) \equiv \frac{1}{2} \cdot s(0, 2) \equiv -q_2,  
\end{align}
\end{subequations}

\noindent
supplying a slightly different proof of some important special cases of (\ref{eq:sPrime}) and (\ref{eq:sDoublePrime}) already stated. The first of these is evaluated mod $p^2$ in Vandiver \cite{Vandiver}, p.\ 112.

\subsection{The case $M=3$}
The case $M=3$ of the first few rows of (\ref{eq:collectArrayTwo}) will likewise be needed below. Let $x = N/M$; then:

\begin{subequations}

\begin{multline}
s(0, 3x) + s(x, 3x) + s(2x, 3x) \\
\equiv 3\left\{ s(0, 3x) + s(1, 3x) + s(2, 3x)\right\} \\
\equiv 3 \cdot s(0, x)
\end{multline}

\begin{multline} \label{eq:M=3}
s(1, 3x) + s(x+1, 3x) + s(2x+1, 3x) \\
\equiv 3\left\{ s(3, 3x) + s(4, 3x) + s(5, 3x)\right\} \\
\equiv 3 \cdot s(1, x).
\end{multline}

\end{subequations}

\noindent

\subsection{The case $M=N/3$}

In the first three rows of (\ref{eq:collectArrayTwo}), let $N$ be divisible by 3, and $M=N/3$. Then 

\subsubsection*{Corollary 4}

\begin{subequations}

\begin{multline} \label{eq:Corollary4}
s(0, N) + s(3, N) + s(6, N) + \ldots + s(N - 3, N) \\
\equiv \frac{N}{3}\left\{ s(0, N) + s(1, N) + s(2, N) + \ldots + s\left(\frac{N}{3} - 1, N\right)\right\} \\
\equiv \frac{N}{3} \cdot s(0, 3) \equiv -\frac{N}{2} \cdot q_3 \qquad
\end{multline}

\begin{multline}
s(1, N) + s(4, N) + s(7, N) + \ldots + s(N - 2, N) \\
\equiv \frac{N}{3}\left\{ s\left(\frac{N}{3}, N\right) + s\left(\frac{N}{3} + 1, N\right) + s\left(\frac{N}{3} + 2, N\right) + \ldots + s\left(\frac{2N}{3} - 1, N\right)\right\} \\
\equiv \frac{N}{3} \cdot s(1, 3) \equiv 0 \qquad
\end{multline}

\begin{multline} 
s(2, N) + s(5, N) + s(8, N) + \ldots + s(N - 1, N) \\
\equiv \frac{N}{3}\left\{ s\left(\frac{2N}{3}, N\right) + s\left(\frac{2N}{3} + 1, N\right) + s\left(\frac{2N}{3} + 2, N\right) + \ldots + s(N - 1, N)\right\} \\
\equiv \frac{N}{3} \cdot s(2, 3)  \equiv \frac{N}{2} \cdot q_3.  
\end{multline}

\end{subequations}

\noindent
Now let $N=p-1$, so that $p \equiv 1 \pmod{3}$; then the sums $s(k, p-1)$ divide $\left\{ 1, p-1\right\}$ into equal pieces of length 1, so that:

\begin{subequations}
\begin{align}
\frac{1}{1} + \frac{1}{4} + \frac{1}{7} + \ldots + \frac{1}{(p - 3)} & \equiv K(1, 3) \equiv \frac{1}{3} \cdot s(2, 3) \equiv \hphantom{-} \frac{1}{2} \cdot q_3 \\
\frac{1}{2} + \frac{1}{5} + \frac{1}{8} + \ldots + \frac{1}{(p - 2)} & \equiv K(2, 3) \equiv \frac{1}{3} \cdot s(1, 3) \equiv \hphantom{-} 0 \\
\frac{1}{3} + \frac{1}{6} + \frac{1}{9} + \ldots + \frac{1}{(p - 1)} & \equiv K(0, 3) \equiv \frac{1}{3} \cdot s(0, 3) \equiv -\frac{1}{2} \cdot q_3.
\end{align}
\end{subequations}

\noindent
These formulae are an improvement upon some of Glaisher (\cite{Glaisher}, p.\ 18), where only the difference between the first and third rows is evaluated. However, Lerch (\cite{Lerch}, p.\ 476) gives a result from which the first row is easily deduced, and Lehmer (\cite{Lehmer1938}, p.\ 356) gives an equivalent result with the modulus $p^2$.

For $p \equiv 2 \pmod 3$, we obtain a parallel set of relations using (\ref{eq:K}):

\begin{subequations}
\begin{align}
\frac{1}{1} + \frac{1}{4} + \frac{1}{7} + \ldots + \frac{1}{(p - 1)} & \equiv K(2, 3) \equiv \frac{1}{3} \cdot s(1, 3) \equiv \hphantom{-} 0 \\
\frac{1}{2} + \frac{1}{5} + \frac{1}{8} + \ldots + \frac{1}{(p - 3)} & \equiv K(1, 3) \equiv \frac{1}{3} \cdot s(2, 3) \equiv \hphantom{-} \frac{1}{2} \cdot q_3 \\
\frac{1}{3} + \frac{1}{6} + \frac{1}{9} + \ldots + \frac{1}{(p - 2)} & \equiv K(0, 3) \equiv \frac{1}{3} \cdot s(0, 3) \equiv -\frac{1}{2} \cdot q_3.
\end{align}
\end{subequations}

\noindent
These likewise improve upon Glaisher (\cite{Glaisher}, p.\ 18), where only the difference between the second and third rows is evaluated. However, Lerch (\cite{Lerch}, p.\ 476) gives a result from which the middle row is easily deduced. Lehmer (\cite{Lehmer1938}, p.\ 356) gives an equivalent result with the modulus $p^2$.

A few additional formulae in a similar vein but involving only terms in the first half of the range $\left\{ 1, p-1\right\}$ follow easily from the foregoing developments. If $p \equiv 1 \pmod 6$,

\begin{subequations}
\begin{align}
\frac{1}{1} + \frac{1}{4} + \frac{1}{7} + \ldots + \frac{1}{(p - 5)/2} & \equiv \frac{1}{3} \cdot s(4, 6) \equiv -\frac{2}{3} \cdot q_2 \\
\frac{1}{2} + \frac{1}{5} + \frac{1}{8} + \ldots + \frac{1}{(p - 3)/2} & \equiv \frac{1}{3} \cdot s(2, 6) \equiv -\frac{2}{3} \cdot q_2 + \frac{1}{2} \cdot q_3 \\
\frac{1}{3} + \frac{1}{6} + \frac{1}{9} + \ldots + \frac{1}{(p - 1)/2} & \equiv \frac{1}{3} \cdot s(0, 6) \equiv -\frac{2}{3} \cdot q_2 - \frac{1}{2} \cdot q_3.
\end{align}
\end{subequations}

\noindent
If $p \equiv 5 \pmod 6$,

\begin{subequations}
\begin{align}
\frac{1}{1} + \frac{1}{4} + \frac{1}{7} + \ldots + \frac{1}{(p - 4)/2} & \equiv \frac{1}{3} \cdot s(2, 6) \equiv -\frac{2}{3} \cdot q_2 + \frac{1}{2} \cdot q_3 \\
\frac{1}{2} + \frac{1}{5} + \frac{1}{8} + \ldots + \frac{1}{(p - 3)/2} & \equiv \frac{1}{3} \cdot s(4, 6) \equiv -\frac{2}{3} \cdot q_2 \\
\frac{1}{3} + \frac{1}{6} + \frac{1}{9} + \ldots + \frac{1}{(p - 2)/2} & \equiv \frac{1}{3} \cdot s(0, 6) \equiv -\frac{2}{3} \cdot q_2 - \frac{1}{2} \cdot q_3.
\end{align}
\end{subequations}

\noindent
Lehmer (\cite{Lehmer1938}, p.\ 358) gives a result equivalent to the middle rows of these last groupings, with the modulus $p^2$. Zhi-Hong Sun (\cite{ZHSun1992}, pt. 1, Corollary 1.12) gives the first two rows of each group.

\section{The remaining classical formulae}

For the sake of completing an inventory of the ``classical'' formulae, we state the following:

If $p \equiv 1 \pmod 4$,

\begin{subequations}
\begin{align}
\frac{1}{1} + \frac{1}{5} + \frac{1}{9} + \ldots + \frac{1}{(p - 4)}  & \equiv K(1, 4) \equiv \frac{1}{4} \cdot s(3, 4) \equiv \hphantom{-}\frac{3}{4} \cdot q_2 \\
\frac{1}{2} + \frac{1}{6} + \frac{1}{10} + \ldots + \frac{1}{(p - 3)} & \equiv K(2, 4) \equiv \frac{1}{4} \cdot s(2, 4) \equiv -\frac{1}{4} \cdot q_2 \\
\frac{1}{3} + \frac{1}{7} + \frac{1}{11} + \ldots + \frac{1}{(p - 2)} & \equiv K(3, 4) \equiv \frac{1}{4} \cdot s(1, 4) \equiv \hphantom{-} \frac{1}{4} \cdot q_2 \\
\frac{1}{4} + \frac{1}{8} + \frac{1}{12} + \ldots + \frac{1}{(p - 1)} & \equiv K(0, 4) \equiv \frac{1}{4} \cdot s(0, 4) \equiv -\frac{3}{4} \cdot q_2,
\end{align}
\end{subequations}

\noindent
Stern (\cite{Stern}, p.\ 187) gives the first and third rows, and Glaisher (\cite{Glaisher}, pp.\ 22--23) gives the first three rows. Lehmer (\cite{Lehmer1938}, p.\ 358) gives a result equivalent to the first row, mod $p^2$.

If $p \equiv 3 \pmod 4$,

\begin{subequations}
\begin{align}
\frac{1}{1} + \frac{1}{5} + \frac{1}{9} + \ldots + \frac{1}{(p - 2)}  & \equiv K(3, 4) \equiv \frac{1}{4} \cdot s(1, 4) \equiv \hphantom{-} \frac{1}{4} \cdot q_2 \\
\frac{1}{2} + \frac{1}{6} + \frac{1}{10} + \ldots + \frac{1}{(p - 1)} & \equiv K(2, 4) \equiv \frac{1}{4} \cdot s(2, 4) \equiv -\frac{1}{4} \cdot q_2 \\
\frac{1}{3} + \frac{1}{7} + \frac{1}{11} + \ldots + \frac{1}{(p - 4)} & \equiv K(1, 4) \equiv \frac{1}{4} \cdot s(3, 4) \equiv \hphantom{-} \frac{3}{4} \cdot q_2 \\
\frac{1}{4} + \frac{1}{8} + \frac{1}{12} + \ldots + \frac{1}{(p - 3)} & \equiv K(0, 4) \equiv \frac{1}{4} \cdot s(0, 4) \equiv -\frac{3}{4} \cdot q_2,
\end{align}
\end{subequations}

\noindent
Stern (\cite{Stern}, p.\ 187) gives the first and third rows, while Glaisher (\cite{Glaisher}, pp.\ 22--23) gives the first row only. Lehmer (\cite{Lehmer1938}, p.\ 358) gives a result equivalent to the third row, mod $p^2$.

If $p \equiv 1 \pmod 6$,

\begin{subequations}
\begin{align}
\frac{1}{1} + \frac{1}{7} + \frac{1}{13} + \ldots + \frac{1}{(p - 6)} & \equiv K(1, 6) \equiv \frac{1}{6} \cdot s(5, 6) \equiv \hphantom{-} \frac{1}{3} \cdot q_2 + \frac{1}{4} \cdot q_3 \\
\frac{1}{2} + \frac{1}{8} + \frac{1}{14} + \ldots + \frac{1}{(p - 5)} & \equiv K(2, 6) \equiv \frac{1}{6} \cdot s(4, 6) \equiv -\frac{1}{3} \cdot q_2 \\
\frac{1}{3} + \frac{1}{9} + \frac{1}{15} + \ldots + \frac{1}{(p - 4)} & \equiv K(3, 6) \equiv \frac{1}{6} \cdot s(3, 6) \equiv \hphantom{-}\frac{1}{3} \cdot q_2 - \frac{1}{4} \cdot q_3 \\
\frac{1}{4} + \frac{1}{10} + \frac{1}{16} + \ldots + \frac{1}{(p - 3)} & \equiv K(4, 6) \equiv \frac{1}{6} \cdot s(2,6) \equiv -\frac{1}{3} \cdot q_2 + \frac{1}{4} \cdot q_3 \\
\frac{1}{5} + \frac{1}{11} + \frac{1}{17} + \ldots + \frac{1}{(p - 2)} & \equiv K(5, 6) \equiv \frac{1}{6} \cdot s(1,6) \equiv \hphantom{-} \frac{1}{3} \cdot q_2 \\
\frac{1}{6} + \frac{1}{12} + \frac{1}{18} + \ldots + \frac{1}{(p - 1)} & \equiv K(0, 6) \equiv \frac{1}{6} \cdot s(0,6) \equiv -\frac{1}{3} \cdot q_2 - \frac{1}{4} \cdot q_3.
\end{align}
\end{subequations}

\noindent
Lehmer (\cite{Lehmer1938}, p.\ 358) gives a result equivalent to the first row, mod $p^2$.

If $p \equiv 5 \pmod 6$,

\begin{subequations}
\begin{align}
\frac{1}{1} + \frac{1}{7} + \frac{1}{13} + \ldots + \frac{1}{(p - 4)} & \equiv K(5, 6) \equiv \frac{1}{6} \cdot s(1, 6) \equiv \hphantom{-} \frac{1}{3} \cdot q_2  \\
\frac{1}{2} + \frac{1}{8} + \frac{1}{14} + \ldots + \frac{1}{(p - 3)} & \equiv K(4, 6) \equiv \frac{1}{6} \cdot s(2, 6) \equiv -\frac{1}{3} \cdot q_2 + \frac{1}{4} \cdot q_3 \\
\frac{1}{3} + \frac{1}{9} + \frac{1}{15} + \ldots + \frac{1}{(p - 2)} & \equiv K(3, 6) \equiv \frac{1}{6} \cdot s(3, 6) \equiv \hphantom{-} \frac{1}{3} \cdot q_2 - \frac{1}{4} \cdot q_3 \\
\frac{1}{4} + \frac{1}{10} + \frac{1}{16} + \ldots + \frac{1}{(p - 1)} & \equiv K(2, 6) \equiv \frac{1}{6} \cdot s(4,6) \equiv -\frac{1}{3} \cdot q_2 \\
\frac{1}{5} + \frac{1}{11} + \frac{1}{17} + \ldots + \frac{1}{(p - 6)} & \equiv K(1, 6) \equiv \frac{1}{6} \cdot s(5,6) \equiv \hphantom{-} \frac{1}{3} \cdot q_2 + \frac{1}{4} \cdot q_3 \\
\frac{1}{6} + \frac{1}{12} + \frac{1}{18} + \ldots + \frac{1}{(p - 5)} & \equiv K(0, 6) \equiv \frac{1}{6} \cdot s(0,6) \equiv -\frac{1}{3} \cdot q_2 - \frac{1}{4} \cdot q_3.
\end{align}
\end{subequations}

\noindent
Lehmer (\cite{Lehmer1938}, p.\ 358) gives a result equivalent to the fifth row, mod $p^2$.

\section{Vandiver's corollary to Lerch's theorem} \label{Vandiver_corollary}

\noindent
Lerch's theorem was rediscovered in 1917 by Vandiver (\cite{Vandiver}, p.\ 114), who gives it in the equivalent form

\begin{equation} \label{eq:Vandiver}
\sum_{k=1}^{N-1} \sum_{j=1}^{\lfloor\frac{kp}{N}\rfloor} \frac{1}{j} := \sum_{j=1}^{N-1} H_{\lfloor jp/N \rfloor} \equiv -N \cdot q_p(N).
\end{equation}

\noindent
which we have seen in section \ref{eq:HarmonicSeries} above. Vandiver points out the following corollary formed by subtracting from this congruence the same congruence for $N+1$:

\begin{equation} \label{eq:VandiverDifference}
\begin{split}
(N & + 1) \cdot q_p(N+1) - N \cdot q_p(N)  \\
& \equiv s\left(N, N\left(N+1\right)\right) \\
& + s\left(2N, N\left(N+1\right)\right) + s\left(2N+1, N\left(N+1\right)\right) \\
& + s\left(3N, N\left(N+1\right)\right) + s\left(3N+1, N\left(N+1\right)\right) + s\left(3N+2, N\left(N+1\right)\right) \\
& \dots \\
& + s\left(N^2, N\left(N+1\right)\right) + \dots + s\left(N^2 + N - 1, N\left(N+1\right)\right) \\
& \equiv s(1, 2) \\
& + 2 \cdot s\left(N, N\left(N+1\right)\right) \\
& + 2 \cdot \left\{ s\left(2N, N\left(N+1\right)\right) + s\left(2N+1, N\left(N+1\right)\right) \right\} \\
& \dots \\
& + 2 \cdot \left\{ s\left(\lfloor N/2 \rfloor \cdot N, N\left(N+1\right)\right) + \dots + s\left(\lfloor N/2 \rfloor \cdot N + \lfloor N/2 \rfloor - 1, N\left(N+1\right)\right) \right\} \\
\end{split}
\end{equation}

\noindent
Although Vandiver does not write out any specific instances of the above, it can be seen that this corollary permits in certain cases the evaluation of $s\left(N, N\left(N+1\right)\right)$ in terms of sums with smaller $N$; thus we have

\begin{multline*}
3 \cdot q_p(3) - 2 \cdot q_p(2) \equiv s(1, 2) + 2 \cdot s(2, 6) \\
\end{multline*}
\begin{multline*}
4 \cdot q_p(4) - 3 \cdot q_p(3) \equiv s(1, 2) + 2 \cdot s(3, 12) \\
\end{multline*}
\begin{multline*}
5 \cdot q_p(5) - 4 \cdot q_p(4) \equiv s(1, 2) + 2 \cdot \left\{ s(4, 20) + s(4, 10) \right\} \\
\end{multline*}
\begin{multline*}
6 \cdot q_p(6) - 5 \cdot q_p(5) \equiv s(1, 2) + 2 \cdot \left\{ s(5, 30) + s(5, 15) \right\} \\
\end{multline*}
\begin{multline*}
7 \cdot q_p(7) - 6 \cdot q_p(6) \equiv s(1, 2) + 2 \cdot \left\{ s(6, 42) + s(6, 21) + s(6, 14) \right\} \\
\end{multline*}
\begin{multline*}
9 \cdot q_p(9) - 8 \cdot q_p(8) \equiv s(1, 2) + 2 \cdot \left\{ s(8, 72) + s(8, 36) + s(8, 24) + s(8, 18) \right\}.
\end{multline*}

\noindent
Full simplification of (\ref{eq:VandiverDifference}) in this fashion requires that $N(N+1)$ be divisible by all the integers up to $\lfloor N/2 \rfloor$, so for example the expression for $8 \cdot q_p(8) - 7 \cdot q_p(7)$ fails to simplify. It will be noticed that the leading term on the right-hand side after $s(1, 2)$ has obvious evaluations such as $s(2, 6) = s(0, 2) - s(0, 3)$, $s(3, 12) = s(0, 3) - s(0, 4)$, and $s(4, 20) = s(0, 5) - s(0, 4)$; or, in general, $s\left(N, N\left(N+1\right)\right) = s(0, N+1) - s(0, N)$. Despite the formal interest of Vandiver's corollary, we doubt that it permits the evaluation of any particular sum that cannot be evaluated in a more obvious way.

\section{Skula's sharpening of Lerch's formula when $N$ is even}

Up to this point, our results may be regarded as a fairly mild generalization of Lerch's formula. In this section, however, we give a new derivation of a transformation which yields significantly improved results. Consider (\ref{eq:Corollary3} \& b), and let: 

\begin{gather*}
a = \textrm{terms in the first row with } k < N/2 \\
b = \textrm{terms in the first row with } k \geq N/2 \\
c = \textrm{terms in the second row with } k < N/2 \\
d = \textrm{terms in the second row with } k \geq N/2.
\end{gather*}

\noindent
Now $a+b$, $c+d$ are already defined in (\ref{eq:Corollary3}), while by definition $a+c \equiv s(0, 2) \equiv -2 \cdot q_2$, $b+d \equiv s(1, 2) \equiv 2 \cdot q_2$, $a \equiv -d$, and $c \equiv -b$. Thus, we have enough information to solve $a$ and $b$ as follows, with $a$ corresponding to row (a) and $b$ to row (b):

\subsubsection*{Theorem 2}

\begin{subequations}

\begin{multline} \label{eq:Theorem2a}
s(0, N) + s(2, N) + s(4, N) + \ldots + s(2 \cdot \lfloor\frac{(N - 1)}{4}\rfloor, N) \\
\equiv \frac{N + 2}{4}\cdot s(0, 2)  \equiv -\frac{N + 2}{2} \cdot q_2 \qquad
\end{multline}

\begin{multline} \label{eq:Theorem2b}
s(1, N) + s(3, N) + s(5, N) + \ldots + s(2 \cdot \lfloor\frac{(N - 3)}{4}\rfloor + 1, N) \\
\equiv -\frac{N - 2}{4}\cdot s(0, 2)  \equiv \frac{N - 2}{2} \cdot q_2.
\end{multline}

\end{subequations}

\noindent
It should be noted that Skula (\cite{Skula2008}, p.\ 8, Corollary 2.4) proved by a somewhat different technique a result equivalent to the second row; and as the sum of (\ref{eq:Theorem2a}) and (\ref{eq:Theorem2b}) is by definition $s(0, 2)$, the value of the first row is an obvious consequence of Skula's result. Nevertheless, we feel that this theorem warrants a closer look because much of its interest lies in the way the results generated by the two rows supplement one another. In the left-hand sides, the values of $s(k, N)$ are simply those with $k$ of the appropriate parity and strictly less than $N/2$. When $N \equiv 0 \pmod{4}$, the number of terms in the left-hand sides of the two rows above is the same; when $N \equiv 2 \pmod{4}$, the number of terms in the left-hand side of the second row is one less than that in the first row. All this will be clearer if we write $k$ in its even form throughout and dovetail the values produced by the two rows:

\begin{equation} \label{eq:SkulaTree}
\resizebox{11cm}{!} {
\begin{tabular} { @{}r @{}c @{}l @{}c @{}r }
~                             & ~        & $s(2,4)$                                   & $\equiv$ &         $-q_2$ \\
                     $s(0,2)$ & $\equiv$ & $s(4,6)$                                   & $\equiv$ & $-2 \cdot q_2$ \\
                     $s(0,4)$ & $\equiv$ & $s(4,8) + s(6,8)$                          & $\equiv$ & $-3 \cdot q_2$ \\
            $s(0,6) + s(2,6)$ & $\equiv$ & $s(6,10) + s(8,10)$                        & $\equiv$ & $-4 \cdot q_2$ \\
            $s(0,8) + s(2,8)$ & $\equiv$ & $s(6,12) + s(8,12) + s(10,12)$             & $\equiv$ & $-5 \cdot q_2$ \\
$s(0,10) + s(2,10) + s(4,10)$ & $\equiv$ & $s(8,14) + s(10,14) + s(12,14)$            & $\equiv$ & $-6 \cdot q_2$ \\
$s(0,12) + s(2,12) + s(4,12)$ & $\equiv$ & $s(8,16) + s(10,16) + s(12,16) + s(14,16)$ & $\equiv$ & $-7 \cdot q_2$ \\
\ldots\ldots\ldots\           & ~        & \ldots\ldots\ldots\                        & ~        &     \ldots     \\
\end{tabular}
}
\end{equation}

\noindent
These conditions are much stronger than those of Lerch (\ref{eq:LerchSplit} or \ref{eq:LerchDiff}), and are particularly parsimonious in the cases where $N$ is oddly even (6, 10, 14, etc.). The case of $N=6$ was solved by Lehmer in 1938 \cite{Lehmer1938}, and that of $N=10$ by Skula in 2008 \cite{Skula2008} (see section \ref{cases_10} below). The cases of $N=8$ and $N=10$ can be proved directly from the evaluations in Zhi-Hong Sun (\cite{ZHSun1992}, pt.\ 3, Theorems 3.1 and 3.3. respectively).

Furthermore, if in (\ref{eq:Theorem2a}) and (\ref{eq:Theorem2b}) we let $N=p-1$; then again the sums $s(k, p-1)$ divide ${0, p-1}$ into equal pieces of length 1, so that: 

\begin{subequations} \label{eq:Glaisher}
\begin{align}
\frac{1}{1} + \frac{1}{3} + \frac{1}{5} + \ldots + \frac{1}{\left\{ 2 \cdot \lfloor\frac{(p - 2)}{4}\rfloor + 1\right\}} & \equiv s^\prime(0, 2) \equiv \frac{1}{4} \cdot s(0, 2) \equiv -\frac{1}{2} \cdot q_2 \\
\frac{1}{2} + \frac{1}{4} + \frac{1}{6} + \ldots + \frac{1}{\left\{ 2 \cdot \lfloor\frac{(p - 4)}{4}\rfloor + 2\right\}} & \equiv s{''}(0, 2) \equiv \frac{3}{4} \cdot s(0, 2) \equiv -\frac{3}{2} \cdot q_2.
\end{align}
\end{subequations}

\noindent
The first row of \ref{eq:Glaisher} supplies a new proof of a theorem of Glaisher (\cite{Glaisher}, p.\ 23, \S40, where one of the versions of the formula is printed with a missing coefficient).

The fact that Theorem 2 results in subsets of $\frac{1}{4}$ of the terms in $\left\{ 1, p-1\right\}$ which are evaluable in terms of Fermat quotients suggests the question of whether subsets of $\frac{1}{6}$ or even of $\frac{1}{12}$ of the terms might be amenable to such treatment, as is the case for some of the individual values of $s(k, N)$ in Table 1. However, attempts to isolate a collection of non-consecutive terms $s(k, N)$ with the $k$ in arithmetic progression and comprising only $\frac{1}{6}$ of the terms in $\left\{1, p-1\right\}$ reveal that the result is evaluable in such a fashion only when the values of $k$ span the entire range. Unlike the formulae leading to the solutions in (\ref{eq:Theorem2a}) and (\ref{eq:Theorem2b}), here the formulae connecting the six pieces leave their values under-determined. Even should additional, as yet undiscovered relations exist, the difficulty can be shown to be insurmountable in general. If we attempt to evaluate the sums of $s(k, N)$ for every third value before or after the midpoint $\frac{(p - 1)}{2}$ we confront the sums $s(0, 12) + s(3, 12)$ and $s(6, 12) + s(9, 12) \equiv -s(5, 12) -s(2, 12)$, and conversely if we attempt to evaluate the sums for every second value within a range of length $\frac{(p - 1)}{3}$ we confront the sums $s(0, 12) + s(2, 12)$, $s(4, 12) + s(6, 12) \equiv -2 \cdot s(0, 12) - 10 \cdot q_2 - \frac{9}{2} \cdot q_3$, and $s(8, 12) + s(10, 12) \equiv -s(3, 12) - s(1, 12)$. None of these results can be expressed solely in terms of Fermat quotients because each entails precisely one value of $s(k, 12)$ which cannot be so expressed (see Table 7). And we already know from the work of Dilcher \& Skula, mentioned in the Introduction, that there are no sums involving less than $\frac{1}{4}$ of the consecutive terms in $\left\{ 1, p-1\right\}$.

Thus, in the ``classical'' theory of the Fermat quotient, apart from the exceptional instances $s(2, 12) \equiv -q_2 + \frac{3}{2} \cdot q_3$, and $s(3, 12) \equiv 3 \cdot q_2 - \frac{3}{2} \cdot q_3$, noted by Vandiver in 1917, it is impossible to devise a sum $s(k, N)$ with any regularity of structure having less than $\frac{1}{12}$ of the terms in $\left\{ 1, p-1\right\}$. Wagstaff (\cite{Wagstaff1978}, p.\ 584), writing about essentially the same problem from another perspective, suggested in 1978: ``Perhaps one could prove that no congruence [of this type] \dots can have fewer than a total of about $p/12$ terms in the sum.'' The question is taken up in his 1987 collaboration with Tanner \cite{TannerWagstaff}, in which they report (p.\ 346): ``we are able to prove a \dots congruence \dots with fewer than about $p/22$ terms in the sums, but it is too complicated to describe here.'' However, as they admit, their miminal sums are somewhat contrived and unrepresentative; and the final section of their paper (pp.\ 347--50) constitutes a sketch of a proof that the normal case ``always has at least about $p/12$ terms.''

Whenever we state categorically that a particular sum of reciprocals cannot be represented by a multiple of a Fermat quotient, it is to understood that at least one prime $p$ is known for which only one of the quantities vanishes mod $p$. And while this does not rule out the possibility that the same sum might be representable as a \textit{combination} of multiples of different Fermat quotients, there are reasons, which we hope to develop in a future paper, for believing that this is unlikely except in the cases covered by Table 1.

\section{Some consequences of our results for particular cases of $N$} \label{cases}

In what follows, we have nothing to add to the results for $N = 1, 2, 3, 4, 6$ (all included in Table 1) as surveyed by such authors as Emma Lehmer \cite{Lehmer1938} and Dilcher \& Skula \cite{DilcherSkula1995}. Rather, we shall elaborate upon the implications of our theorems for certain other values of $N$, especially with reference to the vanishing of $q_2$ or $q_3$, and to the interesting question, considered by Lehmer, of whether they can vanish simultaneously (see section \ref{Lehmersproblem} below). As is well known, the failure of the first case of FLT would require the exponent to satisfying both the congruence $q_2 \equiv 0$ of Wieferich (only known solutions 1093, 3511) and the congruence $q_3 \equiv 0$ of Mirimanoff (only known solutions 11, 1006003). Whether it is possible for the same number to satisfy both congruences remains an open question, but we have been able to sharpen Lehmer's criteria somewhat. 

First, however, we must briefly review progress in the determination of Lerch's sums made since his own paper of 1905. Glaisher (\cite{Glaisher4}, p.\ 273), had already noticed the possibility of evaluating certain sums of powers in which the terms in $\left\{ 1, p-1\right\}$ were split into ranges of 2, 3, 4, or 6 equal parts, in terms of Bernoulli numbers, and Lehmer (\cite{Lehmer1938}, p.\ 352) made the perspicacious observation that these numbers ``can be characterized by the fact that their totient does not exceed two.'' (As previously noted, the computability of $s(2, 12)$ and $s(3, 12)$ as noted by Vandiver depends upon special circumstances.)

The significance of this distinction was revealed in 1991 when results with $\phi(N) = 4$ began to appear. H.\,C.\ Williams (\cite{Williams1991}, p.\ 440) shows that the evaluations of $s(1, 5)$, $s(1, 8)$, and $s(1, 12)$ depend respectively upon the Lucas sequences $U_{p-(\frac{5}{p})}(1, -1)$, $U_{p-(\frac{2}{p})}(2, -1)$, and $U_{p-(\frac{3}{p})}(4, 1)$, the first and second of which correspond to the well-known Fibonacci numbers and Pell numbers; here $(\frac{\hspace{0.5em}}{\hspace{0.5em}})$ is the Jacobi symbol. We shall say more below about the evalution of $s(k, 8)$ and $s(k, 12)$. Sun \& Sun (\cite{Sun+Sun1992}, p.\ 385) evaluate $s(0, 10)$ in terms of $q_2$ and Fibonacci numbers, and Zhi-Hong Sun presents a group of similar formulae involving our $K(r,10)$, some also dependent on $q_5$, from which the values of $s(k, 10)$ for other $k$ can be derived (\cite{ZHSun1992}, Theorem 3.1); and expressions for $s(1,10)$ and $s(2,10)$ involving Lucas numbers (Corollary 1.11). In their joint paper of 1996, Granville \& Sun (\cite{GranvilleSun}, p.\ 119), give a systematic discussion of all cases of $s(0, N)$ where $\phi(N) = 4$, \textit{i.\,e.}, for $N = 5, 8, 10, 12$.

There are also a few results for $\phi(N) > 4$. For $\phi(N) = 6$, in 1992 Zhi-Hong Sun (\cite{ZHSun1992}, Corollary 2.4) evaluated $s^{\ast}(0, 9)$, equivalent to $-s(1, 18)$, in terms of $q_2$ and a complex recurrence relation, while interestingly, the larger ranges $s(k, 9)$ (where $\phi(N)$ also equals 6) have still not been evaluated for any value of $k$. In their 1996 paper, Granville \& Sun (\cite{GranvilleSun}, p.\ 121, evaluate $s(0, 7)$, which also has $\phi(N) = 6$. Incidentally, the remark therein that this is ``the first example with $\phi(N) = 6$'' is true if taken solely with respect to evaluations of Bernoulli polynomials, but if one considers the underlying problem of evaluating partial sums of the Harmonic series, then strictly speaking it would appear to be the second result in this category. The only other $N$ with $\phi(N) = 6$ is 14, and so far as we know there are no results for it.

When $\phi(N) = 8$, the only known result is for $N = 16$ and the case $s(0, 16)$, given in 1993 by Zhi-Hong Sun (\cite{ZHSun1992}, pt. 2, Theorem 2.1); this will be quoted below. The same author (\cite{ZHSun1992}, Theorem 3.2, nos.\ 5 and 6) in effect supplies an evaluation of $s^{\ast}(5, 15)$, and $\phi(15) = 8$, but as we have shown elsewhere \cite{DobsonEisenstein}, there is a special explanation for this case. Zhi-Wei Sun (\cite{ZWSun2006}, p.\ 2216) points out that $s(k, N)$ can be obtained by subtraction from known results for certain values of $k$ when $N = 24, 40, 60$ (\textit{i.\,e.} when $N$ has no prime-power divisors other than 2, 4, 8, 3, or 5). These all involve recurrence sequences such as the Fibonacci numbers, and cannot be evaluated solely in terms of Fermat quotients.

We regret to have to report that Dilcher \& Skula (\cite{DilcherSkula1995}, pp.\ 389--390) are in error when they state that there are zeros of $s(k, N)$ with $p < 2000$ for all values of $N$ from 2 to 46 except 5. There is in fact no such zero for $s(0, 7)$, $s(2, 10)$, $s(3, 10)$, $s(4, 10)$, $s(0, 11)$, $s(2, 11)$, $s(3, 11)$, $s(0, 12)$, $s(3, 12)$, $s(5, 12)$, nor for many other instances of $s(k, N)$ with $N \le 46$ (for some data in the case $k=0$ see Table 10 at the end of this paper). However, these errors have no effect on the paper's main thesis.

Below, we note some calculations of $s(0, N)$ for various $N$ and high values of $p$. These were obtained using methods described in \cite{DobsonMatrixVariation}.

\subsection{$N$ = 8} \label{cases_8}

In the first row of (\ref{eq:M=2}), set $x=4$. Then $s(0, 8) + s(4, 8) \equiv 2\left\{ s(0, 8) + s(1, 8)\right\} \equiv 2 \cdot s(0, 4) \equiv -6 \cdot q_2$ which implies $s(0, 8) + 2 \cdot s(1, 8) + s(3, 8) \equiv 0$. In the first row of (\ref{eq:collectArrayTwo}), set $N=8$ and $M=4$, giving $s(0, 8) + s(2, 8) + s(4, 8) + s(6, 8) \equiv 4\cdot s(0, 2)$. When $q_2 \equiv 0$, a straightforward calculation then gives

\begin{displaymath}
s(0, 8) \equiv -s(1, 8) \equiv -s(2, 8) \equiv s(3, 8),
\end{displaymath}

\noindent
and pairwise, each of these relations is a necessary and sufficient condition for the vanishing of $q_2$. While not without theoretical interest, such conditions do not entail fewer terms than those involving $s(k, 4)$.

As to the actual values of $s(k, 8)$, Williams (\cite{Williams1991}, p.\ 440) evaluates

\begin{equation} \label{eq:WilliamsEight}
U_{p-(\frac{2}{p})}(2, -1)/p \equiv \frac{1}{4} \left\{ s(1, 8) + s(2, 8) \right\},
\end{equation}

\noindent
where $U$ is the Pell sequence 1, 2, 5, 12, 29, \dots (OEIS sequence no.\ A000129), and $\left( \frac{2}{p} \right)$ is a Jacobi symbol. (Here and in what follows, we shall always begin any quoted sequences at $n=1$ rather than $n=0$.) With the application of our Corollary 1, the sum in braces in the right-hand side of (\ref{eq:WilliamsEight}) can be evaluated as $s(0, 2) - s(0, 8) - s(3, 8) \equiv -2 \cdot q_2 + 2s(1, 8)$, allowing the values of $s(k, 8)$ to be obtained for all $k$ (see Table 5). We see that in this case, the combined effect of Lerch's and Williams's results is necessary to resolve the values of $s(k, 12)$ for all $k$; but these values are also readily obtainable from the values of $K(r, 8)$ tabulated in Zhi-Hong Sun (\cite{ZHSun1992}, Theorem 3.3). As every value of $s(k, 8)$ entails a Pell number, one would not expect that they should be expressible as simple multiples of $q_2$; and the fact that they generally cannot be is proven by the following cases where they vanish while $q_2$ does not (the calculations have been extended to $p < 1,250,000,000$, or in the case $k=0$ to $p < 6,691,500,000,000$, without finding any further solution):

\begin{center}
\begin{tabular} { r l }
$s(0, 8)$ & $p$ = 269, 8573, 1300709, 11740973, 241078561 \\
$s(1, 8)$ & $p$ = 29 \\
$s(2, 8)$ & $p$ = 193 \\
$s(3, 8)$ & $p$ = 23, 56993. \\
\end{tabular}
\end{center}

\noindent
We can offer no explanation as to why the distribution of zeroes in the tested range is so strikingly uneven.

Since the vanishing of the sums in the right-hand side of (\ref{eq:WilliamsEight}) is a condition for the failure of the first case of FLT, the same must be true of the Pell quotient. This quotient was checked for $p < 1,000,000,000$ by Elsenhans \& Jahnel \cite{ElsenhansJahnel}, who found the following cases for which it vanishes mod $p$: $p$ = 13, 31, 1546463 (see OEIS A238736). Charles R.\ Greathouse IV, an editor of the OEIS, has informed us in a personal communication that he has extended the calculations for the Pell quotient to $p = 10^{10}$ without finding a further solution.

\subsection{$N$ = 16} \label{cases_16}

In the first row of (\ref{eq:collectArrayTwo}), let $M = 4$. Then

\begin{equation}
s(0, 16) + s(4, 16) + s(8, 16) + s(12, 16) \equiv 4s(0, 4).
\end{equation}

\noindent
In the first two rows of (\ref{eq:M=2}), let $x=8$. Then

\begin{subequations}
\begin{align}
s(0, 16) + s(8, 16)  \equiv 2{s(0, 16) + s(1, 16)} & \equiv 2 \cdot s(0, 8) \\
s(1, 16) + s(9, 16)  \equiv 2{s(2, 16) + s(3, 16)} & \equiv 2 \cdot s(1, 8).
\end{align}
\end{subequations}

\noindent
The first row implies $s(0, 16) + 2 \cdot s(1, 16) - s(8, 16) \equiv 0$. When $q_2 \equiv 0$, $2 \cdot s(0, 8) + 2 \cdot s(1, 8) \equiv 2 \cdot s(0, 4) \equiv 0$, so $s(0, 16) + s(8, 16) + s(1, 16) + s(9, 16) \equiv 0$. Adding these expressions gives $2 \cdot s(0, 16) + 3 \cdot s(1, 16) + s(9, 16) \equiv 0$, furnishing a criterion for the vanishing of $q_2$ which entails only $\frac{3}{16}$ of the terms in $\left\{ 1, p-1\right\}$, a slight improvement on that with $N=8$. However, such an improvement does not continue for higher powers of 2, as the relationship to $s(k, 4)$ becomes too tenuous. Theorem 2 generates congruences in four terms which vanish when $q_2 \equiv 0$. 

In 1993, Zhi-Hong Sun (\cite{ZHSun1992}, pt. 2, Theorem 2.1), showed that

\begin{equation} \label{eq:P12}
s(0, 16) \equiv -4 \cdot q_p(2) - 2 \cdot U_{p-(\frac{2}{p})}(2, -1)/p - 8(S-1)/p,
\end{equation}

\noindent
where $U_{p-(\frac{2}{p})}(2, -1)/p$ is the Pell quotient (OEIS A000129), and

\begin{equation} \label{eq:LerchHarmonicSplit}
S = (-1)^{\lfloor p/16 \rfloor + \lfloor p/8 \rfloor} \times
\begin{dcases}
\left( -C_{n-1} - C_{n}           + C_{n+2} \right) & \text{if $p \equiv \pm 1 \pmod{16}$} \\
\left(  C_{n-1} + C_{n} - C_{n+1}           \right) & \text{if $p \equiv \pm 3 \pmod{16}$} \\
\left(  C_{n-1}         - C_{n+1}           \right) & \text{if $p \equiv \pm 5 \pmod{16}$} \\
\left( -C_{n-1}                             \right) & \text{if $p \equiv \pm 7 \pmod{16}$},
\end{dcases}
\end{equation}

\noindent
with

\begin{subequations}
\begin{gather}
C_0 = C_1 = C_2 = C_4 = 0, C_3 = 1, C_5 = 4, C_6 = -1, C_7 = 14; \\
C_n = 8C_{n-2} - 20C_{n-4} + 16C_{n-6} - 2C_{n-8}.
\end{gather}
\end{subequations}

\subsection{$N$ = 12} \label{cases_12}

We begin by noting that the conditions $s(2, 12) \equiv -q_2 + \frac{3}{2} \cdot q_3 \equiv 0$ and $s(3, 12) \equiv 3 \cdot q_2 - \frac{3}{2} \cdot q_3 \equiv 0$, which Lehmer inexplicably overlooked, provide sharp necessary criteria for the simultaneous vanishing of $q_2$ and $q_3$. Of all such criteria, these sums have the smallest ranges. As in the cases of the conditions $s(0, 6) \equiv -2 \cdot q_2 - \frac{3}{2} \cdot q_3 \equiv 0$ and $s(2, 6) \equiv -2 \cdot q_2 + \frac{3}{2} \cdot q_3 \equiv 0$, these criteria are certainly insufficient individually, as proven by the following cases where the sums vanish although neither $q_2$ nor $q_3$ does (the calculations, which for $s(0, 6)$ extend work of Schwindt \cite{Schwindt}, have been extended at least to $p < 1,600,000$ without finding any further solution): 

\begin{center}
\begin{tabular} { r l }
$s(0, 6)$  & $p$ = 61, 1680023, 7308036881 \\
$s(2, 6)$  & $p$ = 73, 83, 681251 \\
$s(2, 12)$ & $p$ = 179, 619, 17807 \\
$s(3, 12)$ & $p$ = 250829. \\
\end{tabular}
\end{center}

\noindent
Because so much is already known in the case $N=12$, we shall only note further that taking $M=2$ in (\ref{eq:collectArrayTwo}) gives:

\begin{subequations} \label{eq:Twelve}
\begin{align}
s(0, 12) + s(6, 12)  \equiv 2\left\{ s(0, 12) + s(1, 12)\right\} & \equiv 2 \cdot s(0, 6) \equiv -4 \cdot q_2 - 3 \cdot q_3 \\
s(1, 12) + s(7, 12)  \equiv 2\left\{ s(2, 12) + s(3, 12)\right\} & \equiv 2 \cdot s(1, 6) \equiv \hphantom{-} 4 \cdot q_2.
\end{align}
\end{subequations}

\noindent
Here, the more interesting relation is the second one, which gives $s(1, 12) \equiv s(4, 12)$ as another necessary and sufficient condition for the vanishing of $q_2$; thus (truistically) the simultaneous vanishing of $q_2$ and $q_3$ would imply

\begin{displaymath}
s(0, 12) \equiv -s(1, 12) \equiv -s(4, 12) \equiv s(5, 12).
\end{displaymath}

\noindent
Frobenius in his paper of 1914 (\cite{Frobenius}, p.\ 676) gives precisely this condition as a prerequisite for the failure of the first case of FLT, but his proof does not entail the theory of the Fermat quotient.

For $N=12$ and $k < N/2$, only $s(2, 12)$ and $s(3, 12)$ can be evaluated solely in terms of Fermat quotients. As to the other values of $s(k, 12)$, Williams (\cite{Williams1991}, p.\ 440) evaluates

\begin{equation} \label{eq:WilliamsTwelve}
U_{p-(\frac{3}{p})}(4, 1)/p \equiv \frac{1}{6} \left( \frac{3}{p} \right) \left\{ s(1, 12) + s(2, 12) + s(3, 12) + s(4, 12) \right\},
\end{equation}

\noindent
where $U$ is the Lucas sequence 1, 4, 15, 56, \dots (OEIS sequence no.\ A001353), and $\left( \frac{3}{p} \right)$ is a Jacobi symbol. This sequence, which does not seem to have any common name, should not be confused with the standard Lucas numbers (OEIS A000204). The sum in braces in the right-hand side of (\ref{eq:WilliamsTwelve}) is equivalent to $s(0, 2) - s(0, 12) - s(5, 12) \equiv -2 \cdot q_2 + 2 \cdot s(1, 12)$, allowing the values of $s(k, 12)$ to be obtained for all $k$ (see Table 7). As in the case $N$ = 8, we see that the combined effect of Lerch's and Williams's results is necessary to resolve the values of $s(k, 12)$ for all $k$. It may be noted that an evaluation of $s(0, 12)$ can also be recognized with some effort in Granville \& Sun (\cite{GranvilleSun}, p.\ 119). 

Since the vanishing of the sums of the right-hand side of (\ref{eq:WilliamsTwelve}) is a condition for the failure of the first case of FLT, the same must be true of the quotient in the left-hand side. This quotient (see OEIS A238490) was checked for $p < 1,000,000,000$ by Elsenhans \& Jahnel \cite{ElsenhansJahnel}, who found only a single case for which it vanishes mod $p$, namely with $p$ = 103. Continuing the search to $p < 25,000,000,000$, we have found the additional case $p$ = 2297860813.

\subsection{$N$ = 24} \label{cases_24}

As recognized by Zhi-Wei Sun (\cite{ZWSun2006}, p.\ 2216), $s(k, 24)$ can be explicitly evaluated by subtraction from known values when $k = 2, 3, 8, 9$, etc., since 

\begin{subequations}
\begin{align}
s(2, 24) \equiv s(0, 8) - s(0, 12) & \equiv s(1, 12) + s(2, 12) - s(1, 8) \\
s(3, 24) \equiv s(0, 6) - s(0, 8)  & \equiv s(1, 8) - s(2, 12) \\
s(8, 24) \equiv s(2, 6) - s(3, 8)  & \equiv s(2, 8) - s(3, 12) \\
s(9, 24) \equiv s(3, 8) - s(5, 12) & \equiv s(3, 12) + s(4, 12) - s(2, 8).
\end{align}
\end{subequations}

\noindent
See Table 8 at the end for the actual values. In (\ref{eq:M=2}), let $x=12$. Then the relations in the \textit{third} and \textit{fourth} rows are: 

\begin{subequations}
\begin{align}
s(2, 24) + s(14, 24) & \equiv 2\left\{ s(4, 24) + s(5, 24)\right\} \equiv 2 \cdot s(2, 12) \\
s(3, 24) + s(15, 24) & \equiv 2\left\{ s(6, 24) + s(7, 24)\right\} \equiv 2 \cdot s(3, 12),
\end{align}
\end{subequations}

\noindent
which will clearly vanish if it is possible for $q_2$ and $q_3$ to vanish simultaneously. 

The general congruence in the paper by H.\,C.\ Williams (\cite{Williams1991}, p.\ 439, eq.\ 4.6) when his $d$ (our $N$) = 24 takes the form

\begin{equation} \label{eq:24-1}
\begin{split}
s(1, 24) + s(2, 24) & + s(3, 24) + s(4, 24) + 2s(5, 24) + 2s(6, 24) \\
& + s(7, 24) + s(8, 24) + s(9, 24) + s(10, 24) \\
& \equiv 24\thinspace{}\left(\frac{6}{p}\right)\thinspace{}U_{p - (\frac{6}{p})}(10, 1)/p,
\end{split}
\end{equation}

\noindent
where $U(10, 1)$ = 1, 10, 99, 980, 9701, \dots, is OEIS sequence no.\ A004189. We believe the simplest consequence (and with minimal $k$) that can be deduced from this result using Lerch's formula is

\begin{equation} \label{eq:24-2}
s(0, 24) + s(4, 24) \equiv -6q_2 - \frac{3}{2}q_3 -12\thinspace{}\left(\frac{6}{p}\right)\thinspace{}U_{p - (\frac{6}{p})}(10, 1)/p.
\end{equation}

\noindent
Since the vanishing of every quantity in (\ref{eq:24-2}) besides the Lucas quotient is a known condition for the failure of the first case of FLT, then the same must be true for the Lucas quotient. This quotient, with fundamental discriminant $\sqrt{6}$, was tested for $p < 1,000,000,000$ by Elsenhans \& Jahnel \cite{ElsenhansJahnel}, who found only two cases for which it vanishes mod $p$, namely with $p$ = 7 and 523.

In section \ref{Harmonic_numbers} above it was noted that a congruence of Zhi-Wei Sun (\cite{ZWSun2006}, p.\ 2214) supplies evaluations of products over special binomial coefficients with lower index of the form $\lfloor jp/24 \rfloor$, restricted according as $j$ is either a quadratic residue or a nonresidue of 24. Taking the case where $j$ is a residue, and subtracting it from 1, we have

\begin{equation} \label{eq:24-Sun}
\begin{split}
1 - \binom{p-1}{\lfloor p/24 \rfloor} & \binom{p-1}{\lfloor 5p/24 \rfloor} \binom{p-1}{\lfloor 19p/24 \rfloor} \binom{p-1}{\lfloor 23p/24 \rfloor} \\
& \equiv H_{\lfloor p/24 \rfloor} + H_{\lfloor 5p/24 \rfloor} + H_{\lfloor 19p/24 \rfloor} + H_{\lfloor 23p/24 \rfloor} \\
& \equiv -16q_2 - 6q_3 - \left( \frac{6}{p}\right)24\thinspace{}U_{p - (\frac{6}{p})}(10, 1)/p \pmod{p}.
\end{split}
\end{equation}

\noindent
Rewritten in Lerch sums, the sum of the four Harmonic numbers may be expressed as $2 \cdot s(0, 24)$ + $2 \cdot s(4, 24)$ + $2 \cdot s(0, 6)$, so that upon further simplification this congruence proves to be equivalent to that of Williams (\ref{eq:24-2}) above. The case where $j$ is a nonresidue ultimately gives the same result. We have not investigated the extent to which Sun's congruence is independent of Williams's for other cases of $N$. Kuzumaki \& Urbanowicz, in their remarkable paper on binomial coefficients (\cite{KuzumakiUrbanowicz}, p.\ 139), made the first actual advance on (\ref{eq:24-2}), showing in effect that

\begin{equation} \label{eq:KuzumakiUrbanowicz}
\begin{split}
s(0, 24) & \equiv -4q_2 - \frac{3}{2}q_3 - 4\thinspace{}U_{p - (\frac{2}{p})}(2, -1)/p \\
& - 3\thinspace{}\left(\frac{3}{p}\right)\thinspace{}U_{p - (\frac{3}{p})}(4, 1)/p - 6\thinspace{}\left(\frac{6}{p}\right)\thinspace{}U_{p - (\frac{6}{p})}(10, 1)/p.
\end{split}
\end{equation}

\noindent
In conjunction with Williams's result and Lerch's formula, this congruence readily permits the calculation of $s(k, 24)$ for all $k$, and the values are given in Table 8.

\subsection{$N$ = 9} \label{cases_9}

From (\ref{eq:M=3}), $s(0, 9) + s(3, 9) + s(6, 9) \equiv 3\left\{ s(0, 9) + s(1, 9) + s(2, 9)\right\}$ which implies $2 \cdot s(0, 9) + 3 \cdot s(1, 9) + 4 \cdot s(2, 9) - s(3, 9) \equiv 0$, the strongest relation produced by Lerch's theorem other than the ones depending on the relationship with $s(k, 3)$, including those given by \ref{eq:M=3} and other cases of (\ref{eq:collectArrayTwo}) with $n=9$, $M=3$. Although $s(k, 9)$ has not been evaluated for any value of $k$, it is known that in general it cannot be expressed as a simple multiple of $q_3$, as proven by the following cases (apart from the trivial one of $k=4$) where it vanishes while $q_3$ does not (the calculations have been extended at least to $p < 1,638,000$, and in the case $k=0$ to $p < 334,000,000,000$, without finding any further solution):

\begin{center}
\begin{tabular} { r l }
$s(0, 9)$ & $p$ = 677, 6691, 532199813 \\
$s(1, 9)$ & $p$ = 151, 457, 971, 1439, 12613 \\
$s(2, 9)$ & $p$ = 241, 739, 37799 \\
$s(3, 9)$ & $p$ = 97, 58193. \\
\end{tabular}
\end{center}

\subsection{$N$ = 18} \label{cases_18}

We cannot add much to the knowledge of this little-studied case, other than to point out that Corollary 1 gives 

\begin{equation} \label{eq:Eighteen}
s(0, 18) + 2 \cdot s(1, 18) + s(8, 18) \equiv 0, 
\end{equation}

\noindent
while the second row of (\ref{eq:M=3}) gives

\begin{equation}
s(1, 18) + s(7, 18) + s(13, 18) \equiv 6 \cdot q_2,
\end{equation}

\noindent
and the two rows of Theorem 2 give, respectively, 

\begin{subequations}
\begin{align}
s(0, 18) + s(2, 18) + s(4, 18) + s(6, 18) + s(8, 18) & \equiv -10 \cdot q_2 \\
           s(1, 18) + s(3, 18) + s(5, 18) + s(7, 18) & \equiv -8 \cdot q_2. 
\end{align}
\end{subequations}

\noindent
Zhi-Hong Sun (\cite{ZHSun1992}, part 2, Corollary 2.4) in effect evaluates $s(1, 18)$ in terms of a complicated recursion sequence, as may be seen by applying (\ref{eq:sStar_rule}) above to his result. In light of this work, there is no reason to expect that $s(k, 18)$ would vanish with $q_2$ or with $q_3$. Not only does it fail to do so for the two known Wieferich primes and for the two known Mirimanoff primes, but with the possible exceptions of $k = 0, 8$, it cannot in general be expressed as a simple multiple of either of these Fermat quotients, as proven by the following cases where it vanishes while neither of them does (the calculations have been extended at least to $p < 4,000,000$, and in the case $k=0$ to $p < 334,000,000,000$, without finding any further solution):

\begin{center}
\begin{tabular} { r l }
$s(0, 18)$ & $p$ = \ldots \\
$s(1, 18)$ & $p$ = 47, 1777, 217337 \\
$s(2, 18)$ & $p$ = 167 \\
$s(3, 18)$ & $p$ = 1171, 37783 \\
$s(4, 18)$ & $p$ = 137, 251, 1087, 1301, 2111, 5749 \\
$s(5, 18)$ & $p$ = 4177, 1581479 \\
$s(6, 18)$ & $p$ = 108541 \\
$s(7, 18)$ & $p$ = 149, 35267 \\
$s(8, 18)$ & $p$ = \ldots. \\
\end{tabular}
\end{center}

\noindent
The apparent scarcity of zeroes of $s(0, 18)$ and $s(8, 18)$ against those for $s(1, 18)$ has no obvious explanation, as all three figure in the most restrictive relation known (\ref{eq:Eighteen}), and the first two are not known to be more highly constrained than $s(4, 18)$, which has more zeroes than any other value of $s(k, 18)$ in the range tested. Nor is it apparent why the distribution of zeroes in the tested range is so strikingly different from that for $N=8$.

\subsection{$N$ = 5} \label{cases_5}

\noindent
We shall not attempt to treat this case in any detail, for as previously noted, $s(k, 5)$ cannot generally be expressed in terms of Fermat quotients. We merely note that from Lerch's formula (\ref{eq:Lerch}), if $q_5 \equiv 0$, then

\begin{displaymath}
5q_5 \equiv -4s(0, 5) - 2s(1, 5) \equiv 0,
\end{displaymath}

\noindent
which implies that

\begin{displaymath}
-2s(0, 5) - s(1, 5) \equiv 0.
\end{displaymath}

\noindent
From this condition and the obvious relationships between $s(k, 5)$ and $s(k, 10)$, it is straightforward to deduce that if $q_5 \equiv 0$, then

\begin{displaymath}
s(0, 5) \equiv s(4, 10)
\end{displaymath}

\noindent
and

\begin{displaymath}
s(1, 5) \equiv s(1, 10) \equiv -s(3, 10).
\end{displaymath}

\noindent
If we assume the vanishing of $q_5$, then $-2s(0, 5) \equiv s(1, 5)$, so that if one of the sums vanishes, so does the other.

The case $N = 5$ is interesting as the only instance of $s(k, N)$ with $N \le 46$ where no nontrivial zero is known for \textit{any} value of $k$. The evaluations of $s(k, 5)$ are shown in Table 4 below. From these, it is apparent that a zero of $s(1, 5)$ would be a Wall-Sun-Sun prime (sometimes called a Fibonacci-Wieferich prime), i.\,e.\ a prime dividing its Fibonacci quotient $F_{p-(\frac{5}{p})}/p$ (see \cite{Sun+Sun1992}), while a zero of $s(0,5)$ would have $q_5 \equiv -F_{p-(\frac{5}{p})}/p$, where $(\frac{5}{p})$ = 1 or 2 according as $p \equiv \pm 1$ or $\pm 2$ mod 5. It was in fact proven by Vandiver in 1914 \cite{Vandiver1914} that the vanishing of $q_5$ and $s(0, 5)$ were necessary conditions for the first case of FLT, but he did not know the evaluation of $s(0, 5)$ found by Williams and Sun in the 1990s.

\subsection{$N$ = 10} \label{cases_10}

An early result for this value of $N$ was given in terms of the Fibonacci quotient by H.\,C. Williams (\cite{Williams1982}, p.\ 369), who in effect evaluated $s^{\ast}(0, 5)$, equivalent to $-s(1,10)$. Later, the work of Zhi-Hong Sun (\cite{ZHSun1992}, Corollary 1.11 and Theorem 3.1) provided explicit evaluations of $K(r, 10)$, and thus indirectly of $s(k, 10)$, for every value of $k$ (see Table 6). Yet it is nonetheless interesting to consider the relations which pertain among these sums. Indeed Skula (\cite{Skula2008}, pp.\ 9-10) made a special study of this case, and gives $s(0, 10) + 2\cdot s(1, 10) + s(4, 10) \equiv 0$, which corresponds to our Corollary 1 with $x=5$, and $2 \cdot s(0, 10) + 3 \cdot s(1, 10) + 2 \cdot s(2, 10) + 3 \cdot s(3, 10) + 2 \cdot s(4, 10) \equiv 0$, which corresponds to our Theorem 1 with $M=5$. These results can be easily read out of Sun's work (Theorem 3.1) although he does not state them explicitly.

Now Zhi-Hong Sun (\cite{ZHSun1992}, pt. 3, Theorem 3.2, nos. 3 and 5), proved two formulae which together yield the surprising result $s^{\ast}(1, 5) \equiv - s^{\ast}(0, 3)$. This apparently anomalous relationship between sums neither of whose $N$ values divides the other becomes less mysterious when rewritten in the form $s(1, 10) + s(3, 10) \equiv -s(0, 6) - s(2, 6)$, which reveals it as a direct consequence of our Theorem 2. When $q_2 \equiv 0$ the right side of this congruence vanishes, giving $s(1, 10) + s(3, 10) \equiv 0$ as a remarkably compact condition, both necessary and sufficient, for the vanishing of $q_2$. Comparison with Skula's second result then shows that when $q_2 \equiv 0$, we have also $s(0, 10) + s(2, 10) + s(4, 10) \equiv 0$. If in addition we assume that $q_5 \equiv 0$, the relations simplify as follows:

\begin{displaymath}
\begin{tabular} { @{}l @{}c @{}r }
$s(0, 10)$ & $\equiv$ & $3s(0, 5)$ \\
$s(1, 10)$ & $\equiv$ & $s(1, 5)$ \\
$s(2, 10)$ & $\equiv$ & $-4s(0, 5)$ \\
$s(3, 10)$ & $\equiv$ & $-s(1, 5)$ \\
$s(4, 10)$ & $\equiv$ & $s(0, 5)$ \\
\end{tabular}
\end{displaymath}

\noindent
or, more succinctly, 

\begin{displaymath}
4s(0, 10) \equiv -6s(1, 10) \equiv -3s(2, 10) \equiv 6s(3, 10) \equiv 12s(4, 10).
\end{displaymath}

\noindent
The derivations are not difficult, and we skip the details as these relations can be easily inferred from the explicit evaluations of $K(r, 10)$ given in \cite{ZHSun1992}, and which we have adapted for $s(k, 10)$ in Table 6 below.

\subsection{$N$ = 15, 30} \label{cases_15}

\noindent
Note that the sums $s(5, 15)$ and $s(5, 30)$ may be obtained by subtraction. See Table 9 below for details.

\subsection{$N$ = 20} \label{cases_20}

\noindent
Note that (as shown in Table 9 below)

\begin{displaymath}
s(4, 20) = s(0, 4) - s(0, 5) =  -3 \cdot q2 + \frac{5}{4} \cdot q_5 + \frac{5}{4} \cdot F_{p-\epsilon}/p,
\end{displaymath}

\noindent
and

\begin{displaymath}
s(5, 20) = s(2, 10) - s(4, 20) = q_2 - \frac{5}{4} \cdot q_5 + \frac{15}{4} \cdot F_{p-\epsilon}/p,
\end{displaymath}

\noindent
where $F$ is a Fibonacci number and $\epsilon = \left(\frac{5}{p}\right)$. In addition, using Lerch's approach we can obtain relations among sums by setting $M$ = 1, 4, 5, 10 in (\ref{eq:collectArrayTwo}), but no combination of these suffices to isolate any other values of $s(k, 20)$. Setting $M$ = 10, and refering to Table 6 for values of $s(k, 10)$, the second and third rows give

\begin{subequations} \label{eq:20-1}
\begin{align}
s(1, 20) + s(11, 20) \equiv 2s(1, 10) \equiv 2 \cdot q_2 + \frac{5}{2}\thinspace{}F_{p - \epsilon}/p \\
s(2, 20) + s(12, 20) \equiv 2s(2, 10) \equiv -2 \cdot q_2 + 5\thinspace{}F_{p - \epsilon}/p,
\end{align}
\end{subequations}

\noindent
The general congruence in the paper by H.\,C.\ Williams (\cite{Williams1991}, p.\ 439, eq.\ 4.6) when his $d$ (our $N$) = 20 takes the form

\begin{equation} \label{eq:20-2}
s(1, 20) + s(2, 20) - s(7, 20) - s(8, 20) \equiv 10\thinspace{}U_{p - \epsilon}(4, -1)/p,
\end{equation}

\noindent
with $\epsilon$ as before, and where $U(4, -1)$ is OEIS sequence no.\ A001076, but comparison of (\ref{eq:20-2}) with the two rows of (\ref{eq:20-1}) shows that the former is nothing but the sum of the latter, and that therefore

\begin{equation} \label{eq:20-3}
U_{p - \epsilon}(4, -1)/p \equiv \frac{3}{2}\thinspace{}F_{p - \epsilon}/p.
\end{equation}

\noindent
This coincidence is explained by the fact that $U(4, -1)$ (associated with the case $N$ = 20) and the Fibonacci sequence (associated with $N$ = 5) both have $\sqrt{5}$ as their fundamental discriminant (which is just the squarefree part of $N$).

\section{Lehmer's problem} \label{Lehmersproblem}

Although Emma Lehmer was not the first author to pose the question of whether $q_2$ and $q_3$ can vanish simultaneously, her 1938 paper remains the most important contribution to the subject. Indeed, there does not seem to have been much produced since, other that an heuristic argument against the possibility in Lenstra \cite{Lenstra}. However, to the extent that Lehmer develops congruences for Fermat quotients to higher moduli, or derives expressions which cannot be expressed in terms of Lerch's sums, her work is supplemented by the extensive writings of Zhi-Hong Sun, notably by a major recent paper on Bernoulli and Euler numbers \cite{ZHSun2008}.

As previously noted, Lehmer overlooked the conditions involving $s(k, 12)$ discussed above, including the sharpest of all necessary criteria requiring only Fermat quotients, \textit{i.\,e.} $s(2, 12) \equiv -q_2 + \frac{3}{2} \cdot q_3 \equiv 0$ and $s(3, 12) \equiv 3 \cdot q_2 - \frac{3}{2}q_3 \equiv 0$. We have given some comparable conditions involving $N=24$. As to the vanishing of $q_2$ alone, see our (\ref{eq:SkulaTree}), and of $q_3$ alone, our (\ref{eq:Corollary4}); but although each of these results implies an infinite family of conditions, they do not appear to combine in any interesting way.

\section{The effect of the vanishing of $s(k, N)$ on $s(k, 2N)$}

\noindent
Taking (\ref{eq:sStar2}) and (\ref{eq:sStar3}) together, if $s(k, N) \equiv 0$, then

\begin{equation} \label{eq:sStar2Reduced}
s(k, 2N) \equiv s(N-1-k, 2N).
\end{equation}

\noindent
When $N$ is even and $k = \frac{N-2}{2}$, the terms $s(k, 2N) = s\left(\frac{N-2}{2}, 2N\right)$ and $s(N-1-k, 2N) = s\left(\frac{N}{2}, 2N\right)$ are equal, and lie adjacent in the center of the range for $k$ of $\left\{ 0, {N-1}\right\}$. Moreover, when $N$ is only oddly even, these terms have the sum $s\left(\frac{N-2}{2}, N\right)$.

\section{Remark on a result of Dilcher and Skula}

Dilcher and Skula prove in \cite{DilcherSkula1995} that the failure of the first case of Fermat's Last Theorem would imply $s(k, N) \equiv 0$ for all $N \leq 46$ and all $k < N$. For $N$ oddly even, this would give $s\left(\frac{N-2}{2}, 2N\right) + s\left(\frac{N}{2}, 2N\right) \equiv s\left(\frac{N-2}{2}, N\right) \equiv 0$, and since the two terms on the left-hand side are equal and have a vanishing sum, each must itself vanish.

Furthermore, the vanishing of all the sums implied by the proposition of Dilcher and Skula would of course entail the vanishing of $q_2$. On that assumption, adding together (\ref{eq:Theorem2a}) and (\ref{eq:Theorem2b}) and cancelling the vanishing sums $s(0, N) + s(1, N)$, etc., we are left in the case $N \equiv 2 \bmod 4$ with but a single term, $s\left(2 \cdot \lfloor\frac{(N - 1)}{4}\rfloor, N\right) = s\left(\frac{N}{2} - 1, N\right) \equiv 0$, and applying (\ref{eq:M=2}), we find that all the values of $s(k, N)$ may be expressed as multiples of values of $s\left(k, \frac{N}{2}\right)$. In other words, when the condition $s(k, N) \equiv 0$ is proved for an odd $N$ and all $k < N$, the same condition is immediately proved for the case of $2N$. Thus, the result of Dilcher and Skula for $N \leq 46$ automatically extends to all oddly even $N \leq 90$.

This observation complements a result of Cik\'{a}nek \cite{Cikanek}, in which it is shown that the failure of the first case of FLT implies $s(k, N) \equiv 0$ for all $N \leq 94$ and all $k < N$. Cik\'{a}nek's proof requires the additional condition (stated in \S3.4 of the paper) that $p > 5^{(N-1)^{2} (N-2)^{2}/4}$.

It may similarly be noted that when the condition $s(k, N) \equiv 0$ is proved for an odd $N$ and all $k < N$, then by virtue of (\ref{eq:sStar}), it is proved for all values of $s^{\ast}(k, N)$ with the same $N$.

Finally, it should be noted that this argument cannot in general be extended to the case of $2N$ when $N$ is even. As may be seen from an examination of the table for $s(k, 8)$ below (section \ref{cases_8}), its values are affected by the contribution of the Pell numbers even if $s(k, 4)$ vanishes for every value of $k$ (as of course actually happens when $p$ is a Wieferich prime). And as may be seen from the table for $s(k, 12)$, its values for $k = 0, 1, 4, 5$ are affected by the contribution of the Lucas numbers even if $s(k, 6)$ vanishes for every value of $k$.

\clearpage

\begin{table} [hb]
\begin{center}
\caption{Presumably complete list of Lerch's sums (with $k < N/2$) which can be evaluated solely in terms of Fermat quotients}
\label{Table_1}
\begin{tabular}{@{}r | rrr@{}}
$s(0, 1)$  & $ 0          $ & $ ~ $ & $ ~ $ \\
$s(0, 2)$  & $-2 \cdot q_2$ & $ ~ $ & $ ~ $ \\
$s(0, 3)$  & $ ~ $          & $ - $ & $ \frac{3}{2} \cdot q_3$ \\
$s(1, 3)$  & 0 \\
$s(0, 4)$  & $-3 \cdot q_2$ & $ ~ $ & $ ~ $ \\
$s(1, 4)$  & $         q_2$ & $ ~ $ & $ ~ $ \\
$s(0, 6)$  & $-2 \cdot q_2$ & $ - $ & $ \frac{3}{2} \cdot q_3$ \\
$s(1, 6)$  & $ 2 \cdot q_2$ & $ ~ $ & $ ~ $ \\
$s(2, 6)$  & $-2 \cdot q_2$ & $ + $ & $ \frac{3}{2} \cdot q_3$ \\
$s(2, 12)$ & $        -q_2$ & $ + $ & $ \frac{3}{2} \cdot q_3$ \\
$s(3, 12)$ & $ 3 \cdot q_2$ & $ - $ & $ \frac{3}{2} \cdot q_3$ \\
\end{tabular}
\end{center}
\end{table}

\begin{table} [hb]
\begin{center}
\caption{Mod $p^2$ congruences for Lerch's sums, due to Zhi-Hong Sun (\cite{ZHSun2008}, \cite{ZHSun2012}); the $E$ are Euler numbers or polynomials, and $ ( \frac{-1}{\,p} ) = (-1)^{(p-1)/2}$.}
\label{Table_2}
\resizebox{11cm}{!} {
\begin{tabular}{@{}r | rrrrrrrrr@{}}
$s(0, 1)$ &              0 &     &                             &     &                          &     &                              &     & \\
$s(0, 2)$ & $-2 \cdot q_2$ & $+$ & $            p \cdot q_2^2$ &     &                          &     &                              &     & \\
$s(1, 2)$ & $ 2 \cdot q_2$ & $-$ & $            p \cdot q_2^2$ &     &                          &     &                              &     & \\
$s(0, 3)$ &                &     &                             & $-$ & $ \frac{3}{2} \cdot q_3$ & $+$ & $ \frac{3}{4}p \cdot q_3^2 $ & $-$ & $ \frac{1}{9} \left(\frac{p}{3}\right)\, p \cdot E_{p-3}\left(\frac{1}{3}\right) $ \\
$s(1, 3)$ &                &     &                             &     &                          &     &                              &     & $ \frac{1}{3} \left(\frac{p}{3}\right)\, p \cdot E_{p-3}\left(\frac{1}{3}\right) $ \\
$s(2, 3)$ &                &     &                             &     & $ \frac{3}{2} \cdot q_3$ & $-$ & $ \frac{3}{4}p \cdot q_3^2 $ & $-$ & $ \frac{2}{9} \left(\frac{p}{3}\right)\, p \cdot E_{p-3}\left(\frac{1}{3}\right) $ \\
$s(0, 4)$ & $-3 \cdot q_2$ & $+$ & $ \frac{3}{2}p \cdot q_2^2$ &     &                          &     &                              & $-$ & $ \left( \frac{-1}{\,p} \right)\, p \cdot E_{p-3} \hphantom{\left(\frac{1}{3}\right)}$ \\
$s(1, 4)$ & $         q_2$ & $-$ & $ \frac{1}{2}p \cdot q_2^2$ &     &                          &     &                              & $+$ & $ \left( \frac{-1}{\,p} \right)\, p \cdot E_{p-3} \hphantom{\left(\frac{1}{3}\right)}$ \\
$s(2, 4)$ & $        -q_2$ & $+$ & $ \frac{1}{2}p \cdot q_2^2$ &     &                          &     &                              & $+$ & $3\left( \frac{-1}{\,p} \right)\, p \cdot E_{p-3} \hphantom{\left(\frac{1}{3}\right)}$ \\
$s(3, 4)$ & $ 3 \cdot q_2$ & $-$ & $ \frac{3}{2}p \cdot q_2^2$ &     &                          &     &                              & $-$ & $3\left( \frac{-1}{\,p} \right)\, p \cdot E_{p-3} \hphantom{\left(\frac{1}{3}\right)}$ \\
$s(0, 6)$ & $-2 \cdot q_2$ & $+$ & $            p \cdot q_2^2$ & $-$ & $ \frac{3}{2} \cdot q_3$ & $+$ & $ \frac{3}{4}p \cdot q_3^2 $ & $-$ & $  \frac{5}{18} \left(\frac{p}{3}\right)\, p \cdot E_{p-3}\left(\frac{1}{3}\right) $ \\
$s(1, 6)$ & $ 2 \cdot q_2$ & $-$ & $            p \cdot q_2^2$ &     &                          &     &                              & $+$ & $   \frac{1}{6} \left(\frac{p}{3}\right)\, p \cdot E_{p-3}\left(\frac{1}{3}\right) $ \\
$s(2, 6)$ & $-2 \cdot q_2$ & $+$ & $            p \cdot q_2^2$ & $+$ & $ \frac{3}{2} \cdot q_3$ & $-$ & $ \frac{3}{4}p \cdot q_3^2 $ & $+$ & $   \frac{1}{9} \left(\frac{p}{3}\right)\, p \cdot E_{p-3}\left(\frac{1}{3}\right) $ \\
$s(3, 6)$ & $ 2 \cdot q_2$ & $-$ & $            p \cdot q_2^2$ & $-$ & $ \frac{3}{2} \cdot q_3$ & $+$ & $ \frac{3}{4}p \cdot q_3^2 $ & $+$ & $   \frac{2}{9} \left(\frac{p}{3}\right)\, p \cdot E_{p-3}\left(\frac{1}{3}\right) $ \\
$s(4, 6)$ & $-2 \cdot q_2$ & $+$ & $            p \cdot q_2^2$ &     &                          &     &                              & $+$ & $   \frac{7}{6} \left(\frac{p}{3}\right)\, p \cdot E_{p-3}\left(\frac{1}{3}\right) $ \\
$s(5, 6)$ & $ 2 \cdot q_2$ & $-$ & $            p \cdot q_2^2$ & $+$ & $ \frac{3}{2} \cdot q_3$ & $-$ & $ \frac{3}{4}p \cdot q_3^2 $ & $-$ & $ \frac{25}{18} \left(\frac{p}{3}\right)\, p \cdot E_{p-3}\left(\frac{1}{3}\right) $ \\
\end{tabular}
}
\end{center}
\end{table}

\begin{table} [hb]
\begin{center}
\caption{Mod $p^3$ congruences for Lerch's sums, mostly due to Zhi-Hong Sun \cite{ZHSun2008}; the $E$ and $B$ are Euler and Bernoulli numbers, respectively, and $( \frac{-1}{\,p} ) = (-1)^{(p-1)/2}$.}
\label{Table_3}
\resizebox{11cm}{!} {
\begin{tabular}{@{}r | rrrrrrl@{}}
$s(0, 1)$  & $ ~          $ & $ ~ $ & $ ~                       $ & $ ~ $ & $ ~                                           $ & $ - $ & $ \frac{1}{3}p^2 \cdot B_{p-3}$ \\
$s(0, 2)$  & $-2 \cdot q_2$ & $ + $ & $            p \cdot q_2^2$ & $ ~ $ & $ ~                                           $ & $ - $ & $ \;\; p^2 (\frac{2}{3}q_2^3 + \frac{7}{12}B_{p-3})$ \\
$s(1, 2)$  & $ 2 \cdot q_2$ & $ - $ & $            p \cdot q_2^2$ & $ ~ $ & $ ~                                           $ & $ + $ & $ \;\; p^2 (\frac{2}{3}q_2^3 - \frac{1}{4}B_{p-3})$ \\
$s(0, 4)$  & $-3 \cdot q_2$ & $ + $ & $ \frac{3}{2}p \cdot q_2^2$ & $ + $ & $ \medspace \left( \frac{-1}{\,p} \right)\, p \, (E_{2p-4} - 2E_{p-3})$ & $ - $ & $ \;\; p^2(q_2^3 + \frac{7}{12}B_{p-3})$ \\
$s(1, 4)$  & $         q_2$ & $ - $ & $ \frac{1}{2}p \cdot q_2^2$ & $ - $ & $ \medspace \left( \frac{-1}{\,p} \right)\, p \, (E_{2p-4} - 2E_{p-3})$ & $ + $ & $ \frac{1}{3}p^2 \cdot q_2^3$ \\
$s(2, 4)$  & $        -q_2$ & $ + $ & $ \frac{1}{2}p \cdot q_2^2$ & $ - $ & $ 3 \left( \frac{-1}{\,p} \right)\, p \, (E_{2p-4} - 2E_{p-3})$ & $ - $ & $ \frac{1}{3}p^2(q_2^3 + 14B_{p-3})$ \\
$s(3, 4)$  & $ 3 \cdot q_2$ & $ - $ & $ \frac{3}{2}p \cdot q_2^2$ & $ + $ & $ 3 \left( \frac{-1}{\,p} \right)\, p \, (E_{2p-4} - 2E_{p-3})$ & $ + $ & $ \;\; p^2(q_2^3 + \frac{59}{12}B_{p-3})$ \\
\end{tabular}
}
\end{center}
\end{table}

\begin{table} [hb]
\begin{center}
\caption{Values of $s(k, 5)$, where $\mathcal{F}_p = F_{p-(\frac{5}{p})}/p$ is the Fibonacci quotient, and the Jacobi symbol $( \frac{5}{p} )$ = 1 or 2 according as $p \equiv \pm 1$ or $\pm 2$ mod 5. Derived from Williams \cite{Williams1991}, p.\ 440, and Sun \cite{ZHSun1992}, pt.\ 3, Theorems 3.1 and 3.2. The formula for $s(1, 5)$ has also been given by Dilcher \& Skula \cite{DilcherSkula1995}, p.\ 390.}
\label{Table_4}
\begin{tabular}{@{}r | rrrrr@{}}
$s(0, 5)$ & $-\frac{5}{4} \cdot q_5$ & $-$ & $\frac{5}{4} \cdot \mathcal{F}_p$ \\
$s(1, 5)$ & ~                        & ~   & $\frac{5}{2} \cdot \mathcal{F}_p$ \\
\end{tabular}
\end{center}
\end{table}

\begin{table} [hb]
\begin{center}
\caption{Values of $s(k, 8)$, where $\mathcal{P}_p = U_{p-(\frac{2}{p})}(2, -1)/p$ is the Pell quotient (see OEIS sequence no.\ A000129). Derived from Williams \cite{Williams1991}, p.\ 440. Equivalent results also appear in Sun \cite{ZHSun1992}, pt.\ 3, Theorem 3.3.}
\label{Table_5}
\begin{tabular}{@{}r | rrr@{}}
$s(0, 8)$ & $-4 \cdot q_2$ & $-$ & $2 \cdot \mathcal{P}_p$ \\
$s(1, 8)$ & $         q_2$ & $+$ & $2 \cdot \mathcal{P}_p$ \\
$s(2, 8)$ & $        -q_2$ & $+$ & $2 \cdot \mathcal{P}_p$ \\
$s(3, 8)$ & $ 2 \cdot q_2$ & $-$ & $2 \cdot \mathcal{P}_p$ \\
\end{tabular}
\end{center}
\end{table}

\begin{table} [hb]
\begin{center}
\caption{Values of $s(k, 10)$, where $\mathcal{F}_p = F_{p-(\frac{5}{p})}/p$ is the Fibonacci quotient as above. Taken from Sun \cite{ZHSun1992}, pt.\ 3, Theorem 3.1.}
\label{Table_6}
\begin{tabular}{@{}r | rrrrr@{}}
$s(0, 10)$ & $-2 \cdot q_2$ & $-$ & $\frac{5}{4} \cdot q_5$ & $-$ & $\frac{15}{4} \cdot \mathcal{F}_p$ \\
$s(1, 10)$ & $ 2 \cdot q_2$ & ~   &                         & $+$ & $\frac{5}{2}  \cdot \mathcal{F}_p$ \\
$s(2, 10)$ & $-2 \cdot q_2$ & ~   &                         & $+$ & $          5  \cdot \mathcal{F}_p$ \\
$s(3, 10)$ & $ 2 \cdot q_2$ & ~   & ~                       & $-$ & $\frac{5}{2}  \cdot \mathcal{F}_p$ \\
$s(4, 10)$ & $-2 \cdot q_2$ & $+$ & $\frac{5}{4} \cdot q_5$ & $-$ & $\frac{5}{4}  \cdot \mathcal{F}_p$ \\
\end{tabular}
\end{center}
\end{table}

\begin{table} [hb]
\begin{center}
\caption{Values of $s(k, 12)$, where $\mathcal{L}_p = U_{p-(\frac{3}{p})}(4, 1)/p$ is a quotient derived from the Lucas sequence 1, 4, 15, 56, 209, \dots (OEIS sequence no.\ A001353; not to be confused with the standard Lucas numbers), and the Jacobi symbol $( \frac{p}{3} ) = (-1)^{\lfloor p/3 \rfloor}$. Derived from Williams \cite{Williams1991}, p.\ 440.}
\label{Table_7}
\begin{tabular}{@{}r | rrrrr@{}}
$s(0, 12)$ & $-3 \cdot q_2$ & $-$ & $\frac{3}{2} \cdot q_3$ & $-$ & $3 \cdot (\frac{3}{p}) \cdot \mathcal{L}_p$ \\
$s(1, 12)$ & $q_2$          & ~   &                         & $+$ & $3 \cdot (\frac{3}{p}) \cdot \mathcal{L}_p$ \\
$s(2, 12)$ & $-q_2$         & $+$ & $\frac{3}{2} \cdot q_3$ & ~   & ~ \\
$s(3, 12)$ & $3 \cdot q_2$  & $-$ & $\frac{3}{2} \cdot q_3$ & ~   & ~ \\
$s(4, 12)$ & $-3 \cdot q_2$ & ~   & ~                       & $+$ & $3 \cdot (\frac{3}{p}) \cdot \mathcal{L}_p$ \\
$s(5, 12)$ & $q_2$          & $+$ & $\frac{3}{2} \cdot q_3$ & $-$ & $3 \cdot (\frac{3}{p}) \cdot \mathcal{L}_p$ \\
\end{tabular}
\end{center}
\end{table}

\begin{table} [hb]
\begin{center}
\caption{Values of $s(k, 24)$, with $\mathcal{P}_p$ and $\mathcal{L}_p$ as above, and where $U_{p-(\frac{6}{p})}(10, 1)/p$ is a quotient derived from the Lucas sequence 1, 10, 99, 980, 9701, \dots (OEIS sequence no.\ A004189; not to be confused with the standard Lucas numbers). First row shows the component quantities; subsequent rows the corresponding coefficients. Derived from Kuzumaki \& Urbanowicz \cite{KuzumakiUrbanowicz}, p.\ 139.}
\label{Table_8}
\begin{tabular}{@{}r | rrrrr@{}}
~          & $q_2$ & $q_3$ & $\mathcal{P}_p$ &  $(\frac{3}{p}) \mathcal{L}_p$ & $(\frac{6}{p}) U_{p-(\frac{6}{p})}(10, 1)/p$ \\
\hline
$ s(0, 24)$ & $-4$ & $-\frac{3}{2}$ & $-4$ & $-3$ & $-6$ \\
$ s(1, 24)$ & $ 1$ & $           ~$ & $ 4$ & $ ~$ & $ 6$ \\
$ s(2, 24)$ & $-1$ & $ \frac{3}{2}$ & $-2$ & $ 3$ & $ ~$ \\
$ s(3, 24)$ & $ 2$ & $-\frac{3}{2}$ & $ 2$ & $ ~$ & $ ~$ \\
$ s(4, 24)$ & $-2$ & $           ~$ & $ 4$ & $ 3$ & $-6$ \\
$ s(5, 24)$ & $ 1$ & $ \frac{3}{2}$ & $-4$ & $-3$ & $ 6$ \\
$ s(6, 24)$ & $-1$ & $-\frac{3}{2}$ & $-4$ & $ 3$ & $ 6$ \\
$ s(7, 24)$ & $ 4$ & $           ~$ & $ 4$ & $-3$ & $-6$ \\
$ s(8, 24)$ & $-4$ & $ \frac{3}{2}$ & $ 2$ & $ ~$ & $ ~$ \\
$ s(9, 24)$ & $ 1$ & $-\frac{3}{2}$ & $-2$ & $ 3$ & $ ~$ \\
$s(10, 24)$ & $-1$ & $           ~$ & $ 4$ & $-6$ & $ 6$ \\
$s(11, 24)$ & $ 2$ & $ \frac{3}{2}$ & $-4$ & $ 3$ & $-6$ \\
\end{tabular}
\end{center}
\end{table}

\begin{table} [hb]
\begin{center}
\caption{Some additional values of $s(k, N)$ that can be derived from known results by subtraction.}
\label{Table_9}
\resizebox{11cm}{!} {
\begin{tabular}{@{}r | rrrrrrr@{}}
$s(5, 15)$ & ~              & ~   & $\frac{3}{2} \cdot q_3$ & $-$ & $\frac{5}{4} \cdot q_5$ & $+$ & $\frac{5}{4}  \cdot F_{p-(\frac{5}{p})}/p$ \\
$s(4, 20)$ & $-3 \cdot q_2$ & ~   & ~                       & $+$ & $\frac{5}{4} \cdot q_5$ & $+$ & $\frac{5}{4}  \cdot F_{p-(\frac{5}{p})}/p$ \\
$s(5, 20)$ & $         q_2$ & ~   & ~                       & $-$ & $\frac{5}{4} \cdot q_5$ & $+$ & $\frac{15}{4} \cdot F_{p-(\frac{5}{p})}/p$ \\
$s(5, 30)$ & $ 2 \cdot q_2$ & $+$ & $\frac{3}{2} \cdot q_3$ & $-$ & $\frac{5}{4} \cdot q_5$ & $-$ & $\frac{5}{4}  \cdot F_{p-(\frac{5}{p})}/p$ \\
\end{tabular}
}
\end{center}
\end{table}


\begin{table} [hb]
\begin{center}
\caption{Divisors $p$ of Harmonic numbers $H_{\lfloor p/N \rfloor}$ for small $N$. Tests have been carried to $p < 383,950,000$ or higher. Some of these are existing calculations; the divisors $p$ of $H_{\lfloor p/2 \rfloor}$ and $H_{\lfloor p/4 \rfloor}$ are the Wieferich primes (OEIS A001220), while those of $H_{\lfloor p/3 \rfloor}$ are the Mirmimanoff primes (OEIS A014127).}
\label{Table_10}
\begin{tabular}{ r | l || r | l }
$N$ & $p$ & $N$ & $p$ \\
\hline
2  & 1093, 3511                         & 25 & 137                 \\
3  & 11, 1006003                        & 26 & 137, 67939          \\
4  & 1093, 3511                         & 27 & 137, 23669          \\
5  & ---                                & 28 & 20101               \\     
6  & 61, 1680023, 7308036881            & 29 & ---                 \\
7  & 652913                             & 30 & 27089407            \\
8  & 269, 8573, 1300709,                & 31 & ---                 \\
~  & \quad 11740973, 241078561          & 32 & 761                 \\    
9  & 677, 6691, 532199813               & 33 & 761                 \\  
10 & 227, 17539, 4750159                & 34 & 1553                \\ 
11 & 246277, 1156457                    & 35 & 4139, 4481, 4598569 \\
12 & ---                                & 36 & 1297                \\
13 & 43214711, 427794751                & 37 & 1439, 26833         \\
14 & 2267, 6898819                      & 38 & 2473, 3527, 4047089 \\
15 & 134227                             & 39 & 407893              \\
16 & 38723, 38993, 4292543              & 40 & 509, 177553         \\
17 & 590422517                          & 41 & 509, 151883         \\
18 & ---                                & 42 & 509, 190657         \\
19 & 521, 911                           & 43 & ---                 \\
20 & ---                                & 44 & 6967, 27361         \\
21 & 1423, 5693, 5782639, 212084723     & 45 & 609221              \\
22 & 2843                               & 46 & 11731               \\
23 & 137, 264391                        & ~ & ~                    \\
24 & 137, 577, 247421, 307639, 366019,  & ~ & ~                    \\
~  & \quad 5262591617, 31251349243      & ~ & ~                    \\
\end{tabular}
\end{center}
\end{table}

\end{document}